
\ifx\CWIloaded\undefined \let\CWIloaded\relax \else   \fi

\ifx\topglue\undefined \def\topglue{\nointerlineskip\vglue-\topskip\vskip}\fi

\ifx\formloaded\undefined 
  \vsize=24.5 cm \advance \vsize -24pt 
  \hsize=16 cm
  \let\formloaded\relax 
\fi


\catcode`_=11 

\def\let_next{\let\next=}
{\obeylines\gdef\default_address%
 {Universit\'e de Poitiers, D\'epartement de Math\'ematiques,
  UFR Sciences SP2MI, T\'el\'eport~2, BP~30179, %
\ 86962 Futuroscope Chasseneuil Cedex, \ France}}

\newtoks\title \newtoks\MSCcodes \newtoks\CRcodes \newtoks\keywords
\newtoks\titlenote
\newif\if_started

 3

\font\seventeenssbf=cmssbx10 scaled \magstep 3

\font\twelvess=cmss12

\font\tenss=cmss10
\font\tenssit=cmssi10
\font\tenssbf=cmssbx10
\def\sans
{\def\rm{\fam0\tenss}\def\it{\fam\itfam\tenssit}\def\bf{\fam\bffam\tenssbf}\rm
}


\font\ninerm=cmr9 \font\nineit=cmti9 \font\ninebf=cmbx9 \font\ninesl=cmsl9
\font\ninett=cmtt9 \font\ninei=cmmi9 \font\ninesy=cmsy9
\font\niness=cmss9 \font\ninessit=cmssi9 \font\ninessbf=cmssbx10 at 9 pt
\skewchar\ninei='177 \skewchar\ninesy='60

\def\ninept
{\def\rm{\fam0\ninerm}\def\it{\fam\itfam\nineit}\def\bf{\fam\bffam\ninebf}%
 \def\sl{\fam\slfam\ninesl}\def\tt{\fam\ttfam\ninett}\rm
 \def\sans{\def\rm{\fam0\niness}\def\it{\fam\itfam\ninessit}%
           \def\sl{\fam\slfam\ninessit}\def\bf{\fam\bffam\ninessbf}\rm}%
 \textfont0\ninerm \textfont1\ninei \textfont2=\ninesy
 \textfont\itfam\nineit \textfont\slfam\ninesl
 \textfont\bffam\ninebf \textfont\ttfam\ninett
 \normalbaselineskip=11pt \normalbaselines
 \setbox\strutbox=\hbox{\vrule height 8 pt depth 3 pt width 0 pt}%
}


\long\def\centerpar#1{\begingroup
  \advance\leftskip by 0pt plus 1 fil minus 3pt \rightskip=\leftskip
	\parskip=0pt \parindent=0pt \parfillskip=0pt \linepenalty=100
  \noindent\ignorespaces#1\par\endgroup}

\def\good_break{\unskip\penalty-10\ \ignorespaces}

\catcode`@=11
\def\vfootnote#1%
{\insert\footins\bgroup \ninept
 \interlinepenalty=\interfootnotelinepenalty
 \splittopskip=\ht\strutbox 
 \splitmaxdepth=\dp\strutbox 
 \floatingpenalty=20000 
 \leftskip=1pc \rightskip=1pc \spaceskip=0pt \xspaceskip=0pt
 \parindent=0pt \parskip=3pt minus 1pt \parfillskip=0pt plus 1 fil
 \textindent{#1}\footstrut\futurelet\next\f@@t
}
\catcode`@=12

\let\_authors\empty 

%

\def\Author#1{\def\this_author{#1}%
  \ifx\this_author\_eq \let_next\Author 
  \else \let_next\let_next \afterassignment\_author 
  \fi \next}
\def\_author
 {\let\this_addr=\default_address
  \ifx\next\address \let_next\get_address
  \else \ifx\next\email \let_next\get_email
  \else \ifx\next\URL \let_next\get_URL
  \else \add_author
  \fi\fi\fi \next
}

\def\address#1{\errmessage
 {\string\address\space should only follow \string\Author}}
\def\email#1{\errmessage
 {\string\email\space should only follow \string\Author}}
\def\URL#1{\errmessage
 {\string\URL\space should only follow \string\Author}}

\def\get_address{\begingroup\obeylines\get_addr}
\def\get_addr#1%
{\endgroup\def\this_addr{#1}%
  \ifx\this_addr\_cr \let_next\get_address 
  \else\ifx\this_addr\_eq \let_next\get_address 
       \else \let_next\let_next \afterassignment\get_ad
  \fi\fi \next
}
\def\get_ad
{\ifx\next\email \let_next\get_email
 \else \ifx\next\URL \let_next\get_URL
 \else \add_author
 \fi\fi \next
}

\def\get_email#1%
{\def\this_email{#1}%
  \ifx\this_email\_eq \let_next\get_email
  \else \edef\this_addr{\this_addr\endgraf{\tt #1}}%
    \let_next\let_next \afterassignment\after_em
  \fi \next
}
\def\after_em{\ifx\next\URL \let_next\get_URL \else \add_author \fi \next}
\def\get_URL#1%
{\def\this_URL{#1}%
 \ifx\this_URL\_eq \let_next\get_URL
 \else 
  {\def~{\char "7E }%
   \xdef\this_addr{\this_addr\endgraf{\noexpand\tt URL: http://#1}}%
  }\let_next\add_author
 \fi \next
}

\def\_eq{=}{\obeylines\gdef\_cr{^^M}}

\def\add_author 
{\toks0=\expandafter{\_authors}%
 \xdef\_authors
 {\the\toks0 \noexpand\\{\this_author\noexpand\_addr{\this_addr}}}%
}


\def\start_art
{\vfil\eject \ifodd\pageno\else \shipout\vbox to \vsize{}\advancepageno\fi

 \topglue 10mm plus 3 mm minus 2mm
 \centerpar{\let\\=\good_break\seventeenssbf \baselineskip=22pt \the\title}%
 \ifx \_authors\empty \message{Warning: no authors specified.}\else
   \bigskip\centerpar
    {\twelvess
     \def\par{\ifhmode\endgraf\else\medskip\fi}\obeylines
     \def\_addr##1{\smallskip{\ninessit\let\\=\par\noindent##1\endgraf}}%
     \def\\##1{\bigskip##1}\_authors
    }\bigskip
 \fi
 \_startedtrue
}


\outer\def\abstract
{\if_started \errmessage{Abstract must be placed at start of article.}%
 \else \start_art \bigskip \begingroup \ninept \sans \narrower
   \leftline{\hskip\leftskip ABSTRACT}\nobreak\smallskip
	 \everypar{\everypar{}\setbox0=\lastbox}%
 \fi
}

\outer\def\maintext
{\if_started\medskip\endgroup\else\start_art\fi
 \if &\the\MSCcodes\the\CRcodes\the\keywords\the\titlenote&\else
  {\narrower\ninept\sans
   \if &\the\MSCcodes&\else
     \hangindent=4em \noindent
     {\it 1991 Mathematics Subject Classification:\/} \the\MSCcodes.
     \par\fi
   \if &\the\CRcodes&\else
     \hangindent=4em \noindent
     {\it 1991 Computing Reviews Subject Classification:\/} \the\CRcodes.
     \par\fi
   \if &\the\keywords&\else
     \hangindent=2em \noindent{\it Keywords and Phrases:\/} \the\keywords.
     \par\fi
   \if &\the\titlenote&\else
     {\def\\{\unskip\hfil\break\ignorespaces}%
      \hangindent=2em \noindent{\it Note:\/} \the\titlenote.
      \par}\fi
  }\bigskip
 \fi
}

\def\sticksection 
{\vbox\bgroup\parskip=0pt \everypar{\egroup\noindent}}

\def\newreport
{\vfil \eject \secno=0\subsecno=0\proclno=0 \eqnumber=0 \pageno=1 \mark{}
 \title={}\MSCcodes={}\CRcodes={}\keywords={}\titlenote={}
 \let\_authors\empty
 \_startedfalse
}

\catcode`_=8 

\input marxmax

\mathorddef\Part=P2 
\mathorddef\Sym=S\bffam
\def\S#1{\Sym_{#1}}

\mathdef\slid=Rel \triangleright
\mathdef\dils=Rel \triangleleft

\def\lh{<_{\rm h}}
\def\lv{<_{\rm v}}
\let\Card=\#
\def\tr{^{\rm t}}

\mathchardef\edgemap=\delta
\def\edge#1{\edgemap(#1)}
\def\0{\edgemap_0}
\def\dlat{(\Part,{\leq_2})} 

\def\cores#1{{\cal C}_{#1}}
\def\r-{\hbox{$r$-}}

\opdef\diag
\opdef\Dom
\opdef\pos
\opdef\form
\altopdef\hgt{ht}
\opdef\Spin
\opdef\wt
\opdef\LR

\def\ts#1{\tilde\S{#1}}

\def\do#1{d_{\rm o}^#1}
\def\moc#1#2{#1\after#2}

\def\Can#1{{\bf 1}_{#1}}

\def\mathlap#1#2{\rlap{$#1#2$}} 
\def\dublpart#1{#1^{\mathpalette\mathlap\sqcup\sqcap}}

\def\vertr{{\rm T}_0} \def\verconf{{\rm C}_0}
\def\hortr{{\rm T}_1} \def\horconf{{\rm C}_1}

\def\Czbar{\overline C_0}
\newbox\Czbox \setbox\Czbox=\hbox{$\Czbar$}
\setbox0=\hbox{$C$} \ht\Czbox=\ht0 
\dimen0=\wd0 \advance\dimen0 by 1.2pt \wd\Czbox=\dimen0

\def\ar(#1,#2,#3,#4){\({#1\atop#3}\,{#2\atop#4}\)}
\opdef\SST \opdef\SSDT \opdef\Yam \opdef\cq

\newmuskip\bufskip \bufskip=0mu minus 6mu 
\def\[#1]{\mskip\bufskip#1\mskip\bufskip}

\ifx\Tableauxloaded\relax   \else \let\Tableauxloaded\relax \fi


\newcount\cols
{\catcode`,=\active\catcode`|=\active
 \gdef\Young(#1){\hbox{$\vcenter
 {\mathcode`,="8000\mathcode`|="8000
  \def,{\global\advance\cols by 1 &}%
  \def|{\cr
        \multispan{\the\cols}\hrulefill\cr
        &\global\cols=2 }%
  \offinterlineskip\everycr{}\tabskip=0pt
  \dimen0=\ht\strutbox \advance\dimen0 by \dp\strutbox
  \halign
   {\vrule height \ht\strutbox depth \dp\strutbox##
    &&\hbox to \dimen0{\hss$##$\hss}\vrule\cr
    \noalign{\hrule}&\global\cols=2 #1\crcr
    \multispan{\the\cols}\hrulefill\cr%
   }
 }$}}
 \gdef\Skew(#1:#2){\hbox{$\vcenter
 {\mathcode`,="8000\mathcode`|="8000
  \dimen0=\ht\strutbox \advance\dimen0 by \dp\strutbox
  \def\boxbeg{\vbox
    \bgroup\hrule\kern-0.4pt\hbox to\dimen0\bgroup\strut\vrule\hss$}%
  \def\boxend{$\hss\egroup\hrule\egroup}%
  \def,{\boxend\boxbeg}%
  \def|##1:{\boxend\vrule\egroup\nointerlineskip\kern-0.4pt
    \moveright##1\dimen0\hbox\bgroup\boxbeg}%
  \def\\##1\\##2:{\boxend\vrule\egroup\nointerlineskip\kern-0.4pt
    \kern ##1\dimen0\moveright##2\dimen0\hbox\bgroup\boxbeg}%
  \moveright#1\dimen0\hbox\bgroup\boxbeg#2\boxend\vrule\egroup
 }$}}
}


\font\sevenit=cmti7 \font\fiveit=cmti5
\scriptfont\itfam=\sevenit \scriptscriptfont\itfam=\fiveit

\def\smallsquares
{\textfont0=\scriptfont0 \scriptfont0=\scriptscriptfont0
 \textfont1=\scriptfont1 \scriptfont1=\scriptscriptfont1
 \textfont\itfam=\scriptfont\itfam \scriptfont\itfam=\scriptscriptfont\itfam
 \textfont\bffam=\scriptfont\bffam \scriptfont\bffam=\scriptscriptfont\bffam
 \setbox0=\hbox{$($}
 \setbox\strutbox=\hbox{\vrule width 0pt height\ht0 depth\dp0 }
}

\def\bigsquares
{\setbox\strutbox=\hbox
 {\vrule width 0pt height1.3\ht\strutbox depth1.3\dp\strutbox}
}

\newbox\hpair \newbox\vpair
\def\installboxes
{\dimen0=\ht\strutbox \advance\dimen0 by \dp\strutbox
 \setbox0=\vbox{\hrule width \dimen0}
 \setbox1=\rlap{\strut\vrule}
 \setbox\hpair=\rlap{\raise\ht\strutbox\copy0 \lower\dp\strutbox\copy0}
 \setbox\vpair=\rlap{\raise\dimen0\copy1 \kern\dimen0 \copy1}
}

\newcount\ribbonlength
\newcount\joff \newcount\jmax \newcount\imid \newcount\jmid
\def\ribbon #1,#2:#3;#4
{\setbox2=\hbox
 {\count0=0 \count1=0
  \copy1 \rlap{\lower\dp\strutbox\copy0}\if|#4|\else\ribs#4 \fi
  \raise\count1\dimen0\rlap
    {\kern\count0\dimen0 \raise\ht\strutbox\copy0 \copy1\raise0.4pt\copy1}%
  \global\joff=\count0
  \divide\dimen0 by 2 \raise\imid\dimen0\rlap{\kern\jmid\dimen0
    \hbox to 2\dimen0{\hfil$#3$\hfil}}%
 }%
 \lower#1\dimen0\rlap{\kern#2\dimen0\box2}%
 \advance\joff by #2\relax \ifnum \joff>\jmax \global\jmax=\joff \fi
}
\def\ribs#1#2
{\raise\count1\dimen0\rlap
   {\kern\count0\dimen0 \copy\ifnum #1=0 \hpair \else \vpair \fi}%
 \advance\ribbonlength by -2
 \ifnum \ribbonlength=0 \jmid=\count0 \imid=\count1 \fi
 \advance\count#1 by 1
 \ifnum \ribbonlength=1 \jmid=\count0 \imid=\count1 \fi
 \ifnum \ribbonlength>-1 \ifnum \ribbonlength<2
   \advance \jmid by \count0 \advance \imid by \count1
 \fi\fi
 \if |#2|\else \ribs#2 \fi 
}

\def\rtab#1
{\installboxes \ribbonlength=#1 \global\jmax=0
 \setbox0=\hbox{\strut\aftergroup\endrtab \aftergroup}}
\def\endrtab{\advance\dimen0 by \jmax\dimen0 \wd0=\dimen0 \vcenter{\box0}}

\def\arrow#1,#2:#3
{\smash{\lower#1\dimen0\rlap{$\kern#2.5\dimen0
 \if r#3 \,\rightarrow \else
 \if l#3 \llap{$\leftarrow\,$} \else
 \if u#3 \setbox0=\hbox to 0pt{\hss$\uparrow$\hss}\ht0=0pt 
         \raise.5\dimen0 \box0 \else
 \if d#3 \setbox0=\hbox to 0pt{\hss$\downarrow$\hss}\dp0=-2pt
         \lower.5\dimen0 \box0 \else
 \if i#3 \setbox0=\hbox to 0pt{\hss$\nwarrow\,$}\ht0=0pt
         \raise.5\dimen0 \box0 \else
 \if o#3 \setbox0=\hbox to 0pt{$\,\searrow$\hss}\dp0=-2pt
         \lower.5\dimen0 \box0 \else
 \fi\fi\fi\fi\fi\fi
 $}}%
}

\readrefs

\Author={Marc A. A. van Leeuwen}
\email={maavl@mathlabo.univ-poitiers.fr}
\URL{wwwmathlabo.univ-poitiers.fr/~maavl/}

\title={Some bijective correspondences \\ involving domino tableaux}
\MSCcodes={05E10}


\hyphenation{Yamano-uchi hyper-octahedral}

\abstract

We elaborate on the results in \ref{splitting}. We give bijective proofs of a
number of identities that were established there, in particular between the
Yamanouchi domino tableaux and the Littlewood-Richardson tableaux that
correspond to the same tensor product decomposition.

\maintext


\newsection Introduction.

The oldest form of the Robinson-Schensted correspondence, given
in~\ref{Robinson}, associates to any semistandard skew tableaux~$T$ of
shape~$\\/\mu$ and weight~$\alpha$, a pair $(L,P)$ consisting, for some
partition~$\nu$, of a Littlewood-Richardson tableau~$L$ of shape~$\\/\mu$ and
weight~$\nu$, and a semistandard Young tableau~$P$ of shape~$\nu$ and
weight~$\alpha$ (cf.~\ref{L-R intro,~(8)}). In~\ref{splitting}, an
analogue for domino tableaux of this construction was given, associating to a
domino tableau~$D$ a pair consisting of a so-called Yamanouchi domino
tableau~$Y$, and a semistandard Young tableau~$P$. As a result, a
decomposition rule for products of Schur functions similar to the
Littlewood-Richardson rule is derived, counting Yamanouchi domino tableaux
instead of Littlewood-Richardson tableaux; it has the additional advantage
that for the square of a Schur function, its decomposition into its
(representation theoretic) symmetric and alternating parts can be read off.

In this paper we define three bijective constructions complementing these
results. The first associates to a semistandard domino tableau a pair of
ordinary semistandard tableaux called a self-switching tableau pair. This
provides a direct description of the association $D\mapsto P$, of which the
original description in~\ref{splitting}, like that of $T\mapsto P$
in~\ref{Robinson}, is quite indirect, and not even obviously well defined.
The r\^ole of this first construction is therefore analogous to that of
\foreign{jeu de taquin}, which describes the correspondence $T\mapsto P$. Our
second construction is a weight preserving transformation of domino tableaux
of shape~$\\/\mu$ into ordinary tableaux of a shape determined by the
$2$-quotients of $\\$~and~$\mu$, which maps Yamanouchi domino tableaux to
Littlewood-Richardson tableaux; it relates the Littlewood-Richardson rule to
its counterpart using domino tableaux. Our third construction defines a
bijection between certain sets of Yamanouchi domino tableaux, which exhibits
the identity of two different combinatorial expressions for the scalar product
$\<s_\\,\psi^2(s_\mu)>$ given in Corollary~4.3 and Theorem~5.3
of~\ref{splitting}.

To define these bijections, we apply algorithmic constructions defined
elsewhere, notably in \ref{edge sequences} and \ref{L-R intro}. In addition, a
central r\^ole is played by coplactic operations (cf.\ \ref{L-R
intro,~\Sec3}), operating on domino tableaux; an implicit definition of these
operations is contained in the algorithm~\ref{splitting,~7.1}.

This paper is organised as follows. In~\Sec2 we define our main bijective
correspondences, and state (without proof) their fundamental properties.
In~\Sec3 we give the definition of coplactic on domino tableaux used in
this paper, which is in terms of so-called augmented domino tableaux. In~\Sec4
we relate these coplactic operations to those on ordinary tableaux; in
particular the fundamental property of our first bijective correspondence is
established by proving commutation with coplactic operations. In~\Sec5 we do
the same for the other correspondences, which involve moving chains in domino
tableaux.

The current paper depends strongly on notions developed in~\ref{edge
sequences} (of which we only use the special case $r=2$ of dominoes), and
in~\ref{L-R intro}; we shall now recall, and in some cases extend, their
notations. The set of partition is denoted by~$\Part$, and it subset of
$2$-cores by~$\cores2$. The parameters $d^2(\edge{\gamma})$ of a
$2$-core~$\gamma$ always have the form~$(c,-c)$ for~$c\in\Z$, so we shall
denote the $2$-core so parametrised by~$\gamma_c$; explicitly,
$\gamma_c=(k,k-1,\ldots,1)$, where $k=-2c$ if $c\leq0$ and $k=2c-1$ if $c>0$.
We define skew standard tableaux as sequences of partitions (saturated chains
in the Young lattice); unlike in~\ref{edge sequences} we do not specify a set
of entries (an initial subset of~$\N$ is always assumed). We define standard
domino tableaux similarly, using the $2$-rim hook lattice~$\dlat$ instead of
the Young lattice. Semistandard (domino) tableaux are determined by their
standardisation and their weight. We write $\SST(\\/\mu)$ for the set of all
semistandard tableaux of shape~$\\/\mu$, and $\SSDT(\\/\mu)$ for the analogous
set of semistandard domino tableaux; when restricting to entries~$<n$, we
write $\SST(\\/\mu,n)$ and $\SSDT(\\/\mu,n)$. We denote the set of dominoes of
$D\in\SSDT(\\/\mu)$ by~$\Dom(D)$; the set of squares of $T\in\SSDT(\\/\mu)$ is
just the Young diagram~$Y(\\/\mu)$, independently of~$T$. If moreover $s\in
Y(\\/\mu)$ and $x\in\Dom(D)$, we write $T(s)$ respectively $D(x)$ for the
entries of the square~$s$ and domino~$x$, and $\pos(s)$, $\pos(x)$ for their
respective positions (i.e., maximal index of a diagonal meeting the
square/domino). Furthermore we shall use the spin $\Spin(D)$ of~$D$ (half the
number of vertical dominoes in~$\Dom(D)$), and the $2$-sign
$\eps_2(\\/\mu)=(-1)^{2\Spin(D)}$ of its shape.

The affine permutation group $\ts2$ is freely generated by two involutions
$s_0,s_1$; each of them defines a structure of chains on the set of dominoes
of any domino tableau~$D$. Two actions of~$\ts2$ are defined on semistandard
domino tableaux: in the first, denoted by~$\sigma(D)$, the generator~$s_i$
moves all chains of~$D$ for~$s_i$ that are not forbidden (cf.\ \ref{edge
sequences, proposition~4.3.1}); in the second, denoted by~$\moc\sigma D$, it
moves only the subset of open chains (it follows that $\sigma(D)$
and~$\moc\sigma D$ have the same shape). Tableau switching, which is an
involution on pairs~$(S,T)$ of skew semistandard tableaux whose shape fit
together (i.e., $S\in\SST(\mu/\nu)$, $T\in\SST(\\/\mu)$, for some
$\\,\mu,\nu\in\Part$), is denoted by $X(S,T)$, and the relations generated by
inward respectively outward jeu de taquin slides by $T\slid{T'}$ and
$S\dils{S'}$ (when $X(S,T)=(T',S')$, these relations hold; indeed they are
equivalent to having this identity for some~$S,S'$, respectively for
some~$T,T'$).

\beginsection Acknowledgements.

I wish to thank Mark Shimozono for calling my attention to coplactic
operations and their relation to operations on pictures. I also wish to thank
Bernard Leclerc for then many stimulating discussions concerning the
constructions of \ref{splitting}, and in particular their relation with
coplactic graphs.

%
 \par

\newsection Various bijective correspondences for domino tableaux.

In this section we formulate three results about domino tableaux that
complement those found in~\ref{splitting}. We introduce the algorithmic
constructions involved, but postpone proofs of most of their properties. We
use several concepts and constructions that were introduced
in~\ref{splitting}, most notably the concept of Yamanouchi domino tableaux (a
subclass of the semistandard domino tableaux analogous to that of
Littlewood-Richardson tableaux for ordinary semistandard tableaux) and the
bijection of \ref{splitting, Theorem~7.3}. Formally our Yamanouchi domino
tableaux differ from those of \ref{splitting} in that their entries start
from~$0$ (to remain consistent with our other definitions), but this just
amounts to a trivial renumbering.

\subsection Projection from domino tableaux to Young tableaux.

In this subsection we define a weight preserving map from semistandard domino
tableaux to Young tableaux, which will coincide with the projection onto the
second factor after applying the bijection of \ref{splitting, Theorem~7.3}.

As mentioned in the introduction, jeu de taquin can be defined in terms of
tableau switching. On its turn tableau switching is described using so-called
tableau switching families of partitions $(\\^{[i,j]})_{i\in I,j\in J}$, where
$I=\{k,\ldots,l\}$ and~$J=\{m,\ldots,n\}$ are of intervals of~$\Z$. When
$S,T,S',T'$ are skew standard tableaux with $X(S,T)=(T',S')$, then there is
such a family such that $S$~and~$T$ are respectively read off along the left
($j=m$) and bottom ($i=l$) edges of the family, while $T'$~and~$S'$ are
respectively read off along its top ($i=k$) and right ($j=n$) edges. The
conditions for a tableau switching family are such that either $(S,T)$ or
$(T',S')$ determine the entire family. If it happens that $(S,T)=X(S,T)$, we
call $(S,T)$ a self-switching standard tableau pair; the associated tableau
switching family with $I=J$ is symmetric: $\\^{[i,j]}=\\^{[j,i]}$ for all
$i,j\in I$. By the definition of a tableau switching family, such symmetry
arises if and only if the diagonal sub-family $D=(\\^{[i,i]})_{i\in I}$ is a
skew standard domino tableau, and each such domino tableau~$D$ corresponds to
a unique symmetric tableau switching family (the argument is the same as
given in~\ref{eljc, \Sec2.3}). In particular, $D$ determines~$(S,T)$, and
$\pi_0\:D\mapsto(S,T)$ defines a shape preserving bijection~$\pi_0$ from skew
standard domino tableaux to self-switching standard tableau pairs, where the
shape of~$(S,T)$ is taken to be that of the tableau~$S|T$ obtained by joining
the two skew tableaux.

As is the case with tableau switching, the bijection~$\pi_0$ can be extended
to the case of semistandard tableaux; the resulting operation will be denoted
by~$\pi_1$. A pair~$(U,V)$ of semistandard tableaux satisfying $X(U,V)=(U,V)$
is called self-switching; this means that their standardisations form a
self-switching standard tableau pair, and that $\wt U=\wt V$ (which we call
the weight of the pair). For a semistandard domino tableau~$D$ with
standardisation~$D_0$ and weight~$\alpha$, we shall define~$\pi_1(D)$ as the
pair~$(U,V)$ of semistandard tableaux with weight~$\alpha$ and
standardisations~$(S,T)=\pi_0(D_0)$; to justify this, we must verify that
$S$~and~$T$ are compatible with~$\alpha$. Let $d=\\^{[i,i]}-\\^{[i-1,i-1]}$
and~$d'=\\^{[i+1,i+1]}-\\^{[i,i]}$ be successive dominoes of~$D$ with equal
entries, so that $\pos(d)<\pos(d')$ by the definition of semistandard domino
tableaux. This is easily seen to imply that the squares
$s=\\^{[i,i]}-\\^{[i,i-1]}$ and $s'=\\^{[i,i+1]}-\\^{[i,i]}$ satisfy
$\pos(s)<\pos(s')$; therefore the skew standard tableau~$K$ given by row~$i$
of the tableau switching family satisfies compatibility with~$\alpha$ at the
pair of entries under consideration. This compatibility is preserved by jeu de
taquin, so since $S\dils K\dils T$, it holds for $S$ and~$T$ as well. We
have proved:

\proclaim Proposition. \selfswitchprop
The correspondence $\pi_1$ defines a shape and weight preserving bijection
from semistandard domino tableaux to self-switching tableau pairs. \QED

As was done for tableau switching in~\ref{L-R intro, \Sec2.2}, we may deduce a
description of the computation of~$\pi_1(D)$ in terms of sliding entries of
two different colours within the skew diagram. By symmetry, we need only
determine the upper triangular half of the tableau switching family. We label
each vertical difference $\\^{[i+1,j]}-\\^{[i,j]}$ by a red entry, and each
horizontal difference by a blue one; for each colour, the multiset of entries
is that of~$D$. At each stage in the sliding process one has a shuffle of the
weakly increasing sequences of red and blue numbers which, by concatenation of
the associated vertical and horizontal steps, determines a lattice path
through the tableau switching family. Each number in the shuffle corresponds
to an entry in the diagram; among equal numbers of the same colour, the
correspondence preserves left to right order. Whenever a red and blue number
are transposed in the shuffle, the corresponding entries in the diagram are
interchanged if and only if they are in adjacent squares. Initially each
domino of~$D$ is filled with a blue copy of its entry in its inward square and
a red copy in its outward square, corresponding to the shuffle in which each
red number immediately follows the same blue number. This shuffle is
repeatedly modified until all blue numbers precede all red ones, at which
point the blue and red entries define the two component tableaux
of~$\pi_1(D)$. For example, for the computation of
\bigdisplay
\pi_1\left(\,
{\rtab2 \ribbon2,0:1;1 \ribbon2,1:1;1 \ribbon3,0:2;0 \ribbon1,2:2;1
      \ribbon3,2:3;1 \ribbon1,3:3;1 \ribbon4,0:4;0 }
\,\right)=
\left(\,
\vcenter{\offinterlineskip
         \hbox{$\Skew(2:2,3|0:1,1,3|0:2|0:4)$}\hbox{\strut}}
,
\vcenter{\offinterlineskip\hbox{\strut}
\hbox{$\Skew(3:3|1:1,2|1:2,3|0:1,4)$}}
\,\right),
\label(\pioneex)
$$
some intermediate stages, with their shuffles, are (with bold face
representing blue, and italics red):
\bigdisplay
\def\r#1{{\it#1}}\def\b#1{{\bf#1}}
\def\next(#1|#2|#3|#4|#5){\Skew(2:#1|0:#2|0:#3|0:#4|0:#5)}
\everycr{\noalign{\hss}}
\line
{\valign{\hbox{$\strut#$}\tabskip=4pt&\next#\tabskip=0pt\cr
\b1\r1\b1\r1\b2\r2\b2\r2\b3\r3\b3\r3\b4\r4
&(\b2,\b3|\b1,\b1,\r2,\r3|\r1,\r1,\b3|\b2,\r2,\r3|\b4,\r4)\cr
\b1\r1\b1\r1\b2\b2\b3\b3\b4\r2\r2\r3\r3\r4
&(\b2,\b3|\b1,\b1,\b3,\r3|\r1,\r1,\r2|\b2,\r2,\r3|\b4,\r4)\cr
\b1\b1\b2\b2\b3\b3\b4\r1\r1\r2\r2\r3\r3\r4
&(\b2,\b3|\b1,\b1,\b3,\r3|\b2,\r1,\r2|\b4,\r2,\r3|\r1,\r4)\cr
}}
$$
Reordering the shuffle from right to left, as in the example, amounts to
performing inward jeu de taquin slides on the blue entries, into the squares
indicated by the red entries, in decreasing order. Reordering from the left
would amount to performing outward slides on the red entries into the squares
of the blue ones, in increasing order.
Since for $(U,V)=\pi_1(D)$ we have $U\dils V$ by construction, we can define:

\proclaim Definition.
A weight preserving map~$\pi$ from the set of semistandard skew domino
tableaux to the set of semistandard Young tableaux is defined by the condition
$\pi(D)\dils U\dils V$, where $(U,V)=\pi_1(D)$.

\proclaim Theorem. \pithm
For any domino tableau~$D$, the second component of the pair associated to it
by the bijection of \ref{splitting, Theorem~7.3} is equal to~$\pi(D)$. In
particular, if $\wt D=\\\in\Part$, then $D$ is a Yamanouchi domino tableau
if and only if $\pi(D)=\Can\\$, the unique Young tableau with shape and
weight~$\\$.

For instance, for the domino tableau~$D$ of~(\pioneex) we have
$\pi(D)={\smallsquares\smallsquares\Young(1,1,2,3|2,3|4)}$, in agreement with
the Young tableau computed in \ref{splitting, 7.2, example~2}. The method by
which the Young tableau is obtained is entirely different however, and the
proof of theorem~\pithm\ will be a rather indirect one.

\subsection Yamanouchi domino tableaux and Littlewood-Richardson tableaux.

In \ref{splitting, Corollary~4.4} Littlewood-Richardson coefficients are
expressed as the cardinalities of certain sets of Yamanouchi domino tableaux.
We shall exhibit an algorithmic bijection from the corresponding sets of
Littlewood-Richardson tableaux to these sets of domino tableaux. In particular
this defines a partitioning of the set of Littlewood-Richardson tableaux
describing the square of a Schur function, into contributions to the symmetric
and alternating part of the square, by means of the spins of the associated
Yamanouchi domino tableaux. These spins cannot be predicted without performing
the algorithm, which may explain why earlier attempts to describe such a
partitioning by combinatorial means have failed.

We first give a general expression of the multiplication of skew Schur
functions in terms of Yamanouchi domino tableaux (only a special case is
stated explicitly in \ref{splitting}), after introducing some notation for
skew shapes and tableaux characterised by $2$-quotients and $2$-cores.

\proclaim Definition.
For $i=0,1$ let $\\^{(i)}/\mu^{(i)}$ be a skew shape,
$T_i\in\SST(\\^{(i)}/\mu^{(i)})$, and let $\gamma\in\cores2$.
\item{(1)}
   $\cq_2(\gamma,\\^{(0)},\\^{(1)})$ is the unique partition with $2$-core
   $\gamma$ and $2$-quotient $(\\^{(0)},\\^{(1)})$;
\item{(2)}
   $\cq_2(\gamma,\\^{(0)}/\mu^{(0)},\\^{(1)}/\mu^{(1)})=
    \cq_2(\gamma,\\^{(0)},\\^{(1)})/\cq_2(\gamma,\mu^{(0)},\mu^{(1)})$;
\item{(3)} 
   $\cq_2(\gamma,T_0,T_1)$ is the semistandard domino tableau of shape
   $\cq_2(\gamma,\\^{(0)}/\mu^{(0)},\\^{(1)}/\mu^{(1)})$ corresponding to
   $(T_0,T_1)$ under the bijection of \ref{edge sequences, proposition~3.2.2}.

\proclaim Theorem \dueto Carr\'e \& Leclerc. \domprodthm
Let $\chi,\chi'$ be skew shapes and $\gamma\in\cores2$; then
$$
  s_\chi\cdot s_{\chi'}
 =\sum_{\nu\in\Part}\Card{\Yam_2\(\cq_2(\gamma,\chi,\chi'),\nu\)}s_\nu,
$$
where $\Yam_2(\psi,\nu)$ denotes the set of all Yamanouchi domino tableaux of
shape~$\psi$ and weight~$\nu$.

\proof
By \ref{edge sequences, corollary~3.2.3} we have
$$
  \sum_{D\in\SSDT(\cq_2(\gamma,\chi,\chi'),n)}x^{\wt(D)}
 =s_\chi(n)\cdot s_{\chi'}(n),
$$
while by \ref{splitting, Theorem~7.3},
$$
  \sum_{D\in\SSDT(\psi,n)}x^{\wt(D)}
 =\sum_{\nu\in\Part}\left(\Card{\Yam_2(\psi,\nu)}
			  \sum_{T\in\SST(\nu,n)} x^{\wt(T)}\right)
 =\sum_{\nu\in\Part}\Card{\Yam_2(\psi,\nu)} s_\nu(n). \eqno \copy\QEDbox
$$

It follows as a special case that the Littlewood-Richardson coefficient
$c^\nu_{\\,\\'}$ is equal to $\Card{\Yam_2\(\psi,\nu\)}$, where $\psi$ is the
skew shape $\cq_2(\gamma,\\/0,\\'/0)=\cq_2(\gamma,\\,\\')/\gamma$, for an
arbitrary $2$-core $\gamma$. This coefficient is traditionally described as
the cardinality of one of several sets of Littlewood-Richardson tableaux (see
for instance \ref{pictures, 2.6}); we choose the set $\LR(\\*\\',\nu)$ of such
tableaux of the shape~$\\*\\'$ and weight~$\nu$, where $\\*\\'$ is the skew
shape obtained by attaching the diagram of $\\$ to the left and below that
of~$\\'$. The general case of the theorem then has a similar interpretation:
the cardinality of $\Yam_2\(\cq_2(\gamma,\chi,\chi'),\nu\)$ is equal to that
of $\LR(\chi*\chi',\nu)$.
\unskip\footnote\dag
{Strictly speaking, we did not unambiguously define the shape~$\chi*\chi'$,
and if we did, we would not be able to reconstruct $\chi$ and~$ \chi'$
uniquely from it in all cases. Formally we define tableaux
$T_0*T_1\in\SST(\chi*\chi')$ simply as pairs
$(T_0,T_1)\in\SST(\chi)\times\SST(\chi')$, but with the convention that, for
the purpose of reading orders (as in defining $\LR(\chi*\chi',\nu)$), all
entries of~$T_1$ are considered to lie above and to the right of those
of~$T_0$.}
We shall give a bijective proof of this identity: for any $2$-core $\gamma$
and any skew shapes~$\chi,\chi'$ we shall construct a weight preserving
bijection $\LR(\chi*\chi')\to\Yam_2\(\cq_2(\gamma,\chi,\chi')\)$. In fact, the
nature of the construction is such that it simultaneously defines a weight
preserving bijection $\LR(\chi*\chi')\to\Yam_2\(\cq_2(\gamma,\chi',\chi)\)$ as
well; this is remarkable since there is no obvious bijection between
$\LR(\chi*\chi')$ and~$\LR(\chi'*\chi)$. These bijections will be obtained by
restriction of appropriate weight preserving bijections from
$\SST(\chi*\chi')$ to~$\SSDT\(\cq_2(\gamma,\chi,\chi')\)$ and
$\SSDT\(\cq_2(\gamma,\chi',\chi)\)$; to this end we must find such bijections
that map Littlewood-Richardson tableaux precisely to Yamanouchi domino
tableaux.

Without the final condition, such weight preserving bijections are already
provided by the maps $\Sigma_c\:T_0*T_1\mapsto\cq_2(\gamma_c,T_0,T_1)$ and
$\Sigma'_c\:T_0*T_1\mapsto\cq_2(\gamma_c,T_1,T_0)$, for~$c\in\Z$; these do not
however in general map Littlewood-Richardson tableaux to Yamanouchi domino
tableaux. This can be understood from the definitions of such tableaux, which
require the word formed by reading the entries of the tableau in a particular
order to be a lattice permutation: for Littlewood-Richardson tableaux this can
be any valid reading order as described in~\ref{L-R intro, \Sec1.5}, while for
Yamanouchi domino tableaux this is the reverse of the column reading order
of~\ref{splitting} (a Yamanouchi word is a reverse lattice permutation). Each
entry of~$\Sigma_c(T_0*T_1)$ or $\Sigma'_c(T_0*T_1)$ has a matching entry
in~$T_0*T_1$, but the reading order used in the domino tableau does not always
correspond to a valid reading order in~$T_0*T_1$. This even fails for the
components $T_0$ and~$T_1$ individually, but more importantly, their entries
are interleaved in the reading of~the domino tableau, while all entries
of~$T_1$ precede those of ~$T_0$ first factor in any valid reading
of~$T_0*T_1$.

However, if $|c|$ is sufficiently large, the situation is different. According
to \ref{edge sequences, proposition~3.1.2}, a domino $d$ in
$\cq_2(\gamma_c,U,V)$ corresponding to a square $s$ of~$U$ satisfies
$\pos(d)=2(\pos(s)+c)$, while it satisfies $\pos(d)=2(\pos(s)-c)+1$ if it
corresponds to a square $s$ of~$V$. Therefore $\Sigma_c$ preserves the
relative order of positions when $c\ll0$, and $\Sigma'_c$ does so when
$c\gg0$. Here is a concrete exemple to illustrate what happens: we display
$T_0*T_1$, followed by $\Sigma'_3(T_0*T_1)$, \ $\Sigma_{-2}(T_0*T_1)$, \
$\Sigma'_2(T_0*T_1)$, and $\Sigma_{-1}(T_0*T_1)$.
$${\Skew(3:0,0,0|3:1|1:1,2|0:0,3|0:2)}
\qquad\qquad
\def\tab{\smallsquares \rtab2 \lower9\dimen0\rlap{\strut}}
\def\+#1{\llap{$\Sigma_{#1}$\quad}\qquad}
\def\-#1{\llap{$\Sigma'_{#1}$\quad}\qquad}
{\tab
 \ribbon10,0:2;1 \ribbon8,0:0;1 \ribbon7,1:3;1
 \ribbon5,1:1;1 \ribbon4,2:2;1
 \ribbon1,4:1;0 \ribbon0,5:0;0 \ribbon0,7:0;0 \ribbon0,9:0;0
}\-3
{\tab 
 \ribbon9,0:2;1 \ribbon7,0:0;1 \ribbon6,1:3;1
 \ribbon4,1:1;1 \ribbon3,2:2;1
 \ribbon1,3:1;0 \ribbon0,4:0;0 \ribbon0,6:0;0 \ribbon0,8:0;0
}\+{-2}
{\tab 
 \ribbon8,0:2;1 \ribbon6,0:0;1 \ribbon5,1:3;1
 \ribbon3,1:1;1 \ribbon2,2:2;0
 \ribbon1,2:1;0 \ribbon0,3:0;0 \ribbon0,5:0;0 \ribbon0,7:0;0
}\-2
{\tab
 \ribbon7,0:2;1 \ribbon5,0:0;1 \ribbon4,1:3;1
 \ribbon2,1:1;1 \ribbon2,2:1;1 \ribbon2,3:2;1
 \ribbon0,2:0;0 \ribbon0,4:0;0 \ribbon0,6:0;0
}\+{-1}\unskip
$$
In the first two domino tableaux displayed, the dominoes corresponding to
squares of~$T_0$ and~$T_1$ form disjoint subtableaux, the former consisting
entirely of vertical dominoes and the latter of horizontal ones; in both
components the correspondence between squares and dominoes is linear, and
independent of the other component. Then the reading order in~$T_0*T_1$
induced by the column reading order for the domino tableau is a valid one, so
the fact that $T_0*T_1$ is a Littlewood-Richardson tableau implies that the
indicated domino tableaux are Yamanouchi. For the third domino tableaux above
most of these statements loose their validity, and the fourth one is in fact
no longer a Yamanouchi domino tableau. The condition that characterises the
simpler situation in the first two domino tableaux can be stated in terms of
the dominoes $d$~and~$d'$ corresponding respectively to the top-right
square~$s$ of~$T_0$ and the bottom-left square~$t$ of~$T_1$: it is
$\pos(d)\leq\pos(d')-3$ (so that the diagonal with index $\pos(d)+1$ separates
the two components), which in view of the expression given above becomes
$2c\leq\pos(s)-\pos(t)-1$ in case of~$\Sigma_c$, and $2c\geq\pos(t)-\pos(s)+2$
for~$\Sigma'_c$. Indeed, in the example $\pos(s)-\pos(t)=-1-2=-3$, so these
conditions are met for $c\leq-2$ respectively for $c\geq3$. Concluding, our
reasoning leads to the following definition and proposition.

\proclaim Definition.
Let $T_0\in\SST(\\/\mu)$ and $T_1\in\SST(\\'/\mu')$; then $\Sigma_c(T_0*T_1)$
(respectively $\Sigma'_c(T_0*T_1)$) is called a segregated tableau
for~$T_0*T_1$ when $2c\leq-n+1$ (respectively when $2c\geq n$), where
$n=\\_0+(\\')\tr_0$.

\proclaim Proposition. \segrprop
If a domino tableau $\Sigma_c(T_0*T_1)$ or $\Sigma'_c(T_0*T_1)$ is segregated,
then it is a Yamanouchi domino tableau if and only if $T_0*T_1$ is an
Littlewood-Richardson tableau.
\QED

Our goal is now to extend the proposition by replacing $\Sigma_c(T_0*T_1)$
and~$\Sigma'_c(T_0*T_1)$, in case they are not segregated, by appropriate
other domino tableaux of the same shape and weight. To this end observe that
the collection of all domino tableaux $\Sigma_c(T_0*T_1)$
and~$\Sigma'_c(T_0*T_1)$ for $c\in\Z$ forms an orbit for the action
$(\sigma,D)\mapsto\sigma(D)$ of the group~$\ts2$ on semistandard domino
tableaux. Indeed $s_i(\cq_2(\gamma,T_0,T_1))=\cq_2(s_i(\gamma),T_1,T_0)$ for
$i=0,1$ by the description of \ref{edge sequences, proposition~4.3.2}, while
$s_0(\gamma_c)=\gamma_{1-c}$ and $s_1(\gamma_c)=\gamma_{-c}$; the orbit can be
depicted as follows (omitting the arguments~$T_0*T_1$):
$$\def\s#1{\buildrel \textstyle s_#1\over\longleftrightarrow}
\cdots
\Sigma'_2\s0\Sigma_{-1}\s1
\Sigma'_1\s0\Sigma_0\s1
\Sigma'_0\s0\Sigma_1\s1
\Sigma'_{-1}\s0\Sigma_2
\cdots
$$
If we go far enough to the left in this diagram, the domino tableaux
$\Sigma_c$ and~$\Sigma'_c$ will be segregated. We shall choose such a
segregated domino tableau~$S$ in this orbit, and then replace the orbit by the
orbit of~$S$ for the other action $(\sigma,D)\mapsto\moc\sigma D$ of~$\ts2$
(cf.\ \ref{edge sequences, proposition~4.5.2}). Recall that whereas in the
former action application of~$s_i$ amounts to moving all non-forbidden chains
for~$s_i$, this is limited in the latter action to moving the subset of open
chains. Since it is easily seen that a segregated tableau contains only open
chains (both for~$s_0$ and for~$s_1$), the part of the two orbits consisting
of segregated tableaux will be identical, whence our construction is
independent of the choice of~$S$.

\proclaim Definition. \Phidef
For any skew shapes $\chi,\chi'$ and $c\in\Z$, maps $\Phi_c$~and~$\Phi'_c$
from $\SST(\chi*\chi')$ to respectively $\SSDT(\cq_2(\gamma\c,\chi,\chi'))$
and $\SSDT(\cq_2(\gamma\c,\chi',\chi))$ are defined as follows. Let $S$ be a
segregated domino tableau for $T_0*T_1\in\SST(\chi*\chi')$; then
$\Phi_c(T_0*T_1)$ and $\Phi'_c(T_0*T_1)$ are the unique elements in the orbit
of~$S$ for the action $(\sigma,D)\mapsto\moc\sigma D$, of respective shapes
$\cq_2(\gamma,\chi,\chi')$ and $\cq_2(\gamma,\chi',\chi)$: if
$\sigma,\sigma'\in\ts2$ are such that $\sigma(S)=\Sigma_c(T_0*T_1)$ and
$\sigma'(S)=\Sigma'_c(T_0*T_1)$, then $\Phi_c(T_0*T_1)=\sigma\after S$ and
$\Phi'_c(T_0*T_1)=\sigma'\after S$.

As an example we compute $\Phi_1(T_0*T_1)$ for the tableau $T_0*T_1$ of the
previous example. We start with the third domino tableau desplayed earlier,
which is the first non-segregated one, and is part of the orbit for either of
the actions; after this point there is divergence, and the action
$(\sigma,D)\mapsto\moc\sigma D$ proceeds:
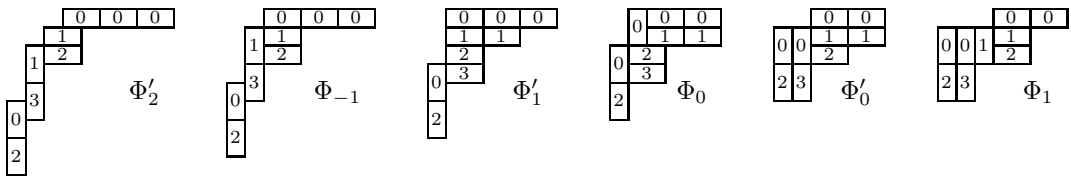
\bigdisplay
\def\tab{\smallsquares \rtab2 \lower8\dimen0\rlap{\strut}}
\def\+#1{\llap{$\Phi_{#1}$\enspace}\qquad}
\def\-#1{\llap{$\Phi'_{#1}$\enspace}\qquad}
{\tab 
 \ribbon8,0:2;1 \ribbon6,0:0;1 \ribbon5,1:3;1
 \ribbon3,1:1;1 \ribbon2,2:2;0
 \ribbon1,2:1;0 \ribbon0,3:0;0 \ribbon0,5:0;0 \ribbon0,7:0;0
}\-2
{\tab
 \ribbon7,0:2;1 \ribbon5,0:0;1 \ribbon4,1:3;1
 \ribbon2,1:1;1 \ribbon2,2:2;0 \ribbon1,2:1;0
 \ribbon0,2:0;0 \ribbon0,4:0;0 \ribbon0,6:0;0
}\+{-1}
{\tab
 \ribbon6,0:2;1 \ribbon4,0:0;1 \ribbon3,1:3;0
 \ribbon1,1:1;0 \ribbon1,3:1;0 \ribbon2,1:2;0
 \ribbon0,1:0;0 \ribbon0,3:0;0 \ribbon0,5:0;0
}\-1
{\tab
 \ribbon5,0:2;1 \ribbon3,0:0;1 \ribbon3,1:3;0
 \ribbon1,2:1;0 \ribbon1,4:1;0 \ribbon2,1:2;0
 \ribbon1,1:0;1 \ribbon0,2:0;0 \ribbon0,4:0;0
}\+0
{\tab
 \ribbon4,0:2;1 \ribbon2,0:0;1 \ribbon4,1:3;1
 \ribbon1,2:1;0 \ribbon1,4:1;0 \ribbon2,2:2;0
 \ribbon2,1:0;1 \ribbon0,2:0;0 \ribbon0,4:0;0
}\-0
{\tab
 \ribbon4,0:2;1 \ribbon2,0:0;1 \ribbon4,1:3;1
 \ribbon2,2:1;1 \ribbon1,3:1;0 \ribbon2,3:2;0
 \ribbon2,1:0;1 \ribbon0,3:0;0 \ribbon0,5:0;0
}\+1\unskip
$$
Note that these are all Yamanouchi domino tableaux. As we shall prove, this is
no coincidence:

\proclaim Theorem. \blowcorethm
By restriction, $\Phi_c$~and~$\Phi'_c$ define, for any~$c\in\Z$, bijections
from $\LR(\chi*\chi')$ respectively to $\Yam_2(\cq_2(\gamma_c,\chi,\chi'))$
and to $\Yam_2(\cq_2(\gamma_c,\chi',\chi))$.

We shall in fact prove the stronger statement that moving any open chain is
compatible with the bijection of \ref{splitting, Theorem~7.3}: on the first
factor (the Yamanouchi domino tableau) the corresponding open chain is moved,
while the second factor (the Young tableau) is unchanged.

Concerning definition~\Phidef, we note the following. While, as indicated in
the example, it is clear where in the orbit each $\Phi_c(T_0*T_1)$ and
$\Phi'_c(T_0*T_1)$ are to be found, we have initially defined these tableaux
as being uniquely determined within the orbit by their shape. When the shapes
$\chi$~and~$\chi'$ of $T_0$~and~$T_1$ are distinct, this is true because the
action of~$\ts2$ on the shapes in the orbit is free; on the other hand when
$\chi=\chi'$, each shape has a stabiliser of order~$2$, since
$\cq_2(\gamma_0,\chi,\chi)$ is stabilised by~$s_1$. In the latter case however
the tableau $\Phi_0(T_0*T_1)$ is also stabilised by~$s_1$ (i.e.,
$\Phi_0(T_0*T_1)=\Phi'_0(T_0*T_1)$), since the fact that the shape does not
change implies that all chains are closed; then the whole orbit is symmetric
and each tableau is unique for its shape. This also implies that
$\Phi_c$~and~$\Phi'_c$ coincide for $\chi=\chi'$.

Another point that can be observed is that the image of the map~$\Phi_c$ for
the shape~$\chi*\chi'$ coincides with the image of~$\Phi'_c$ for~$\chi'*\chi$.
Therefore we may compose the former map with the inverse of the latter, so as
to obtain a bijection $\SST(\chi*\chi')\to\SST(\chi'*\chi)$ which, by
theorem~\blowcorethm, restricts to a bijection
$X_{\Dom}\:\LR(\chi*\chi')\to\LR(\chi'*\chi)$; the following description shows
that $X_{\Dom}(T_0*T_1)$ does not depend on~$c$. One forms segregated tableau
for~$T_0*T_1$, and moves though the non-segregated part of its orbit for the
action $(\sigma,D)\mapsto\moc\sigma D$, to the other (right) end, where the
domino tableau becomes segregated again, but for some other
tableau~$T'_1*T'_0=X_{\Dom}(T_0*T_1)$ of shape~$\chi'*\chi$. This bijection
differs from the ``traditional'' bijection
$X_{\LR}\:\LR(\chi*\chi')\to\LR(\chi'*\chi)$ that is described in detail
in~\ref{path jeu de taquin}, and can be characterised in the language
of~\ref{L-R intro} as tableau switching on companion tableaux. In fact the two
bijections have rather different characteristics: in $X_{\Dom}$ the
shapes~$\chi,\chi'$ play a crucial r\^ole (e.g., $X_{\Dom}$ is the identity
for $\chi=\chi'$, by the remarks above), whereas $X_{\LR}$ only uses some
reading of the tableaux (a lattice permutation); also $X_{\LR}$ does not
naturally extend to a bijection~$\SST(\chi*\chi')\to\SST(\chi'*\chi)$.

\subsection Matching two expressions for $\<s_\\,\psi^2(s_\mu)>$.

Our third construction establishes a bijection corresponding to the identity
$$
  \eps_2(\\)\Card{\Yam_2(\\,\mu)}
 =\sum_{M\in\Yam_2(\dublpart\mu,\\)}(-1)^{|\mu|-\Spin(M)},
 \label(\dublpartid)
$$
where $\dublpart\mu=\cq_2(\emptyset,\mu,\mu)$ is the partition obtained from
the Young diagram of~$\mu$ by scaling up by a factor~$2$ both horizontally and
vertically. Since $\Spin(D)$ denotes half the number of vertical dominoes
in~$D$, and $\eps_2(\\)=\eps_2(\\/\emptyset)=(-1)^{2\Spin(D)}$ for
any~$D\in\SSDT(\\)$, the Yamanouchi domino tableaux in the first member are
counted with a fixed sign determined by the parity of their number of vertical
dominoes, while those in the second member are counted with a varying sign,
determined by the parity of \emph{half} the number of \emph{horizontal}
dominoes (since $\eps_2(\dublpart\mu)=1$, and these tableaux have $2|\mu|$
dominoes altogether).

We recall that both members of~(\dublpartid) are combinatorial expressions for
the number $\<s_\\,\psi^2(s_\mu)>$, where $\psi^2$ is the plethysm operator
that replaces each monomial~$x^\alpha$ by~$x^{2\alpha}$. We review briefly the
derivation of this identity, as it is spread across many sections
of~\ref{splitting}, and since similar arguments will be used in our bijective
proof. Denoting by~$\phi^2$ the dual of the linear operator~$\psi^2$, so that
$\<s_\\,\psi^2(s_\mu)>=\<\phi^2(s_\\),s_\mu>$ for all~$\\,\mu\in\Part$, it is
a classical fact that $\phi^2(s_\chi)=\eps_2(\chi)s_{\chi'}s_{\chi''}$ when
$\chi=\cq_2(\gamma,\chi',\chi'')$ for skew shapes $\chi',\chi''$
(see~\ref{Littlewood modular}; $\phi^2(s_\chi)=0$ if the shape~$\chi$ admits
no domino tableaux). Therefore by theorem~\domprodthm:
$$
   \phi^2(s_\chi)=\eps_2(\chi)\sum_{\nu\in\Part}\Card{\Yam_2(\chi,\nu)}s_\nu
   \label(\phieq)
$$
(cf. \ref{splitting, Corollary~4.3, (5)}), whence the first member
of~(\dublpartid) equals~$\<\phi^2(s_\\),s_\mu>$. The second member is obtained
by evaluating in two ways the sum of~$(-1)^{\Spin(D)}x^{\wt(D)}$ as $D$
ranges over $\SSDT(\dublpart\mu,n)$. On one hand, since the map from domino
tableaux to Yamanouchi domino tableaux in \ref{splitting, Theorem~7.3}
preserves the spin, this sum decomposes as
$$ 
  \sum_{\nu\in\Part}
  \left(
    \sum_{\scriptstyle M\in\Yam_2(\dublpart\mu,\nu)
          \atop
	  \scriptstyle T\in\SST(\nu,n)
	 }
    (-1)^{\Spin(M)} x^{\wt(T)}
  \right)=
  \sum_{\nu\in\Part}
  \biggl(\sum_{M\in\Yam_2(\dublpart\mu,\nu)}(-1)^{\Spin(M)}\biggr)s_\nu(n).
  \label(\spingfdecompeq)
$$
On the other hand, one can cancel from the sum all contributions of domino
tableaux~$D$ among whose chains for~$s_1$ (which are all closed) there is
least one that can be moved (cf.\ \ref{edge sequences, proposition~4.3.1}):
the tableau obtained by moving one such chain contributes with an opposite
sign by~\ref{edge sequences, proposition~4.4.1}. What remains are those domino
tableaux $D\in\SSDT(\dublpart\mu,n)$ for which every $2\times2$ block
corresponding to a square of~$\mu$ is occupied by a pair of vertical dominoes
with equal entries, forming a forbidden chain for~$s_1$. These tableaux are in
bijection with ordinary semistandard tableaux~$D'$ of shape~$\mu$; since
$\Spin(D)=|\mu|$ and $\wt D=2\wt D'$, the summation becomes
$$ 
  \sum_{D'\in\SST(\mu,n)}(-1)^{|\mu|}x^{2\wt D'}= 
  (-1)^{|\mu|}\psi^2(s_\mu)(n).
  \label(\spingfcanceleq)
$$
Taking the coefficient of~$s_\\(n)$ in (\spingfdecompeq) and
(\spingfcanceleq), one finds that $\<s_\\,\psi^2(s_\mu)>$ equals the second
member of~(\dublpartid) (cf.~\ref{splitting, Theorem~5.3}), which establishes
that identity since $\<s_\\,\psi^2(s_\mu)>=\<\phi^2(s_\\),s_\mu>$.

Our combinatorial construction corresponding to~(\dublpartid) will consist of
a bijection between tableaux $L\in\Yam_2(\\,\mu)$ and tableaux~$M$ in a
subset~$B_{\\,\mu}$ of $\Yam_2(\dublpart\mu,\\)$, such that~$L$ has half as
many vertical dominoes as $M$ has horizontal dominoes (so that
$\eps_2(\\)=(-1)^{2\Spin(L)}=(-1)^{|\mu|-\Spin(M)}$), together with a proof
that $\sum_{M\in C_{\\,\mu}}(-1)^{\Spin(M)}=0$, where $C_{\\,\mu}$ is the
complement of~$B_{\\,\mu}$ in $\Yam_2(\dublpart\mu,\\)$. That proof is similar
to the argument leading to~(\spingfcanceleq): we define $C_{\\,\mu}$ to be the
subset of~$\Yam_2(\dublpart\mu,\\)$ of tableaux that contain at least one
chain for~$s_1$ that can be moved while preserving the Yamanouchi property.
The contributions of these tableaux to the sum will cancel out, provided we
can show that the following holds:

\proclaim Lemma. \closchindlem
Let $M$ be a Yamanouchi domino tableau, $s\in\{s_0,s_1\}$, and let $S$ be the
set of closed chains~$C$ in~$M$ for~$s$, for which the tableau obtained
from~$M$ by moving~$C$ is again a Yamanouchi domino tableau. Then the tableau
obtained from~$M$ by simultaneously moving the chains of any subset of~$S$ is
also a Yamanouchi domino tableau.

In other words, the set~$S$ itself is not affected by moving any of its
chains. This means that the relation between Yamanouchi domino tableaux of
being obtainable from one another by moving a subset of the chains of~$S$ is
an equivalence relation, whose equivalence classes have size~$2^{\Card{S}}$;
since moving any one chain changes the parity of~$\Spin(M)$, the sum of
$(-1)^{\Spin(M)}$ over such a class is~$0$ if $S\neq\emptyset$. A similar
fact, but with ``Yamanouchi'' replaced by ``semistandard'', was needed in the
derivation of~(\spingfcanceleq), but that fact is simpler: it follows directly
from \ref{edge sequences, proposition~4.3.1}. We note that a remark is made in
\ref{splitting} (after its Lemma~8.5) that appears to claim the validity our
lemma~\closchindlem; however this remark is neither justified nor used there,
and so we shall provide a proof of the lemma below.

Given $L\in\Yam_2(\\,\mu)$, we construct a filling~$M$ of~$Y(\dublpart\mu)$
using an augmentation of domino tableaux similar to that of ordinary tableaux
in \ref{L-R intro,~(15)}, by attaching subscripts called ordinates to the
entries of~$L$. We do this in such a way that if the dominoes $d$ with fixed
entry~$L(d)=i$ are listed by \emph{decreasing} value of~$\pos(d)$ (i.e., from
right to left), then their ordinates increase by unit steps, starting at~$0$.
This being done, the domino containing $i_j$ (entry~$i$ with ordinate~$j$)
determines two dominoes of~$M$, which occupy the $2\times2$ block
$\{2i,2i+1\}\times\{2j,2j+1\}$ of squares in~$Y(\dublpart\mu)$. These two
dominoes will be horizontal if $d$ is vertical, and vice versa; their two
entries are the row numbers of the two squares forming~$d$ (which are equal if
the dominoes are vertical, and different in they are horizontal, in which case
of course the top domino gets the smaller entry). While it is not obvious that
$M$ is a Yamanouchi (or even a semistandard) domino tableau, it \emph{is}
clear that it has twice as many horizontal dominoes as $T$ has vertical
dominoes.

\proclaim Theorem. \yamyamthm
For any $\\,\mu\in\Part$ and $L\in\Yam_2(\\,\mu)$, the filling~$M$
of~$Y(\dublpart\mu)$ constructed above lies in the subset $B_{\\,\mu}$ of\/
$\Yam_2(\dublpart\mu,\\)$ (i.e., none of its chains for~$s_1$ can be moved
without destroying the Yamanouchi property); moreover the construction defines
a bijection $\Yam_2(\\,\mu)\to B_{\\,\mu}$.

As an example, consider $\\=(6,5,3,3,3)$ and $\mu=(4,3,2,1)$. Now there
is just one $L\in\Yam_2(\\,\mu)$, displayed here together with the
corresponding element $M\in\Yam_2(\dublpart\mu,\\)$, both with augmentation:
\bigdisplay
L=\,
{\rtab2
 \ribbon1,0:0_3;1 \ribbon1,1:0_2;1 \ribbon0,2:0_1;0 \ribbon0,4:0_0;0
 \ribbon2,0:1_2;0 \ribbon2,2:1_1;1 \ribbon1,3:1_0;0
 \ribbon4,0:2_1;1 \ribbon3,1:2_0;0
 \ribbon4,1:3_0;0
}\qquad
M=\,
{\rtab2
 \ribbon1,6:1_0;0 \ribbon0,6:0_0;0  \ribbon1,4:1_1;0 \ribbon0,4:0_1;0
 \ribbon1,3:0_2;1 \ribbon1,2:0_3;1  \ribbon1,1:0_4;1 \ribbon1,0:0_5;1
 \ribbon3,4:2_0;1 \ribbon3,5:2_1;1  \ribbon3,2:2_2;0 \ribbon2,2:1_2;0
 \ribbon3,1:1_3;1 \ribbon3,0:1_4;1
 \ribbon5,2:4_0;0 \ribbon4,2:3_0;0  \ribbon5,1:3_1;1 \ribbon5,0:3_2;1
 \ribbon7,1:4_1;1 \ribbon7,0:4_2;1
}
$$
As can be seen, in addition to the fact that the entries of each pair of
dominoes in a $2\times2$ block of~$M$ give the row numbers of the squares of
the corresponding domino of~$L$, their ordinates give the column numbers. The
four tableaux of the remainder~$C_{\\,\mu}$ of $\Yam_2(\dublpart\mu,\\)$ form
a single equivalence class, for which the set~$S$ of lemma~\closchindlem\
for~$s_1$ consists of two closed chains, each occupying one $2\times2$
block:
\bigdisplay {\smallsquares\installboxes \global\dimen7=3\dimen0 }
\def\spin#1{\lower\dimen7\llap{$\Spin=#1$}}
\matrix{
 {\rtab2
  \ribbon1,6:1_0;0 \ribbon0,6:0_0;0  \ribbon1,4:1_1;0 \ribbon0,4:0_1;0
  \ribbon1,3:0_2;1 \ribbon1,2:0_3;1  \ribbon1,1:0_4;1 \ribbon1,0:0_5;1
  \ribbon3,4:3_0;0 \ribbon2,4:2_0;0  \ribbon3,2:2_1;0 \ribbon2,2:1_2;0
  \ribbon3,1:1_3;1 \ribbon3,0:1_4;1
  \ribbon5,2:4_0;0 \ribbon4,2:3_1;0  \ribbon5,0:3_2;0 \ribbon4,0:2_2;0
  \ribbon7,1:4_1;1 \ribbon7,0:4_2;1
 }
 \spin4
& \longleftrightarrow &
 {\rtab2
  \ribbon1,6:1_0;0 \ribbon0,6:0_0;0  \ribbon1,4:1_1;0 \ribbon0,4:0_1;0
  \ribbon1,3:0_2;1 \ribbon1,2:0_3;1  \ribbon1,1:0_4;1 \ribbon1,0:0_5;1
  \ribbon3,4:3_0;0 \ribbon2,4:2_0;0  \ribbon3,2:1_2;1 \ribbon3,3:2_1;1
  \ribbon3,1:1_3;1 \ribbon3,0:1_4;1
  \ribbon5,2:4_0;0 \ribbon4,2:3_1;0  \ribbon5,0:3_2;0 \ribbon4,0:2_2;0
  \ribbon7,1:4_1;1 \ribbon7,0:4_2;1
 }
 \spin5
\cr
\Big\updownarrow & & \Big\updownarrow \cr \noalign{\nobreak\smallskip}
 {\rtab2
  \ribbon1,6:1_0;0 \ribbon0,6:0_0;0  \ribbon1,4:1_1;0 \ribbon0,4:0_1;0
  \ribbon1,3:0_2;1 \ribbon1,2:0_3;1  \ribbon1,1:0_4;1 \ribbon1,0:0_5;1
  \ribbon3,4:3_0;0 \ribbon2,4:2_0;0  \ribbon3,2:2_1;0 \ribbon2,2:1_2;0
  \ribbon3,1:1_3;1 \ribbon3,0:1_4;1
  \ribbon5,3:4_0;1 \ribbon5,2:3_1;1  \ribbon5,0:3_2;0 \ribbon4,0:2_2;0
  \ribbon7,1:4_1;1 \ribbon7,0:4_2;1
 }
 \spin5
& \longleftrightarrow &
 {\rtab2
  \ribbon1,6:1_0;0 \ribbon0,6:0_0;0  \ribbon1,4:1_1;0 \ribbon0,4:0_1;0
  \ribbon1,3:0_2;1 \ribbon1,2:0_3;1  \ribbon1,1:0_4;1 \ribbon1,0:0_5;1
  \ribbon3,4:3_0;0 \ribbon2,4:2_0;0  \ribbon3,2:1_2;1 \ribbon3,3:2_1;1
  \ribbon3,1:1_3;1 \ribbon3,0:1_4;1
  \ribbon5,3:4_0;1 \ribbon5,2:3_1;1  \ribbon5,0:3_2;0 \ribbon4,0:2_2;0
  \ribbon7,1:4_1;1 \ribbon7,0:4_2;1
 }
 \spin6
\cr
}
$$
We have added ordinates in these tableaux as well; one can see that the blocks
of the chains in~$S$ have pairs $\{1_2,2_1\}$ and $\{3_1,4_0\}$ of
entry-ordinate combinations that do not correspond to adjacent squares.

 \par

\newsection Coplactic operations.

We shall now study in detail the algorithm defining the bijection of
\ref{splitting, Theorem~7.3}, which is given in~\ref{splitting, 7.1}. It is
modelled after Robinson's algorithm, which is defined in \ref{Robinson,
\Sec5}, but is more clearly described in \ref{Macdonald, I~9}. The basic
steps of Robinson's algorithm are the coplactic operations, for a the detailed
discussion of which we refer to~\ref{L-R intro}. Similar operations can be
defined for domino tableaux as well; these are the basic steps of the first
algorithm mentioned above, and will also be called coplactic operations. We
recall that the basic definition of coplactic operations is the one for words,
which is then transferred to ordinary semistandard tableaux by applying it to
words obtained by reading the entries of the tableaux in some specific order.
The transfer to semistandard domino tableaux is more involved, for two related
reasons: there is no obvious way to apply the valid reading orders for
ordinary tableaux to domino tableaux, and coplactic operations may involve
rearrangement of dominoes (sometimes a coplactic operation is applicable, but
no proper result can be obtained by just changing the entry of a domino).

The method of definition adopted in~\ref{splitting} is to use a particular
reading order called the column reading of a domino tableau, which detects
whether a coplactic operation can be applied, and if so makes a change to one
entry as a first approximation of the operation; in case the result is not a
semistandard domino tableau, a transformation called~$R_1$ or a sequence of
transformations called~$R_2$ is applied, yielding the desired result. This
method is not without difficulties; the effect of the transformations on
further coplactic operations is hard to perceive, as they entail permutations
of the letters of the column reading. In particular it is not obvious that the
``coplactic graph'' associated to a domino tableau is always isomorphic to
that of some ordinary tableau, or equivalently, of some word; in any case one
cannot always take the column reading word. Yet this property is essential for
establishing the bijection of~\ref{splitting, Theorem~7.3}, since the second
component of the pair it associates to a domino tableau~$D$ can be described
as ``the Young tableau that is in the same place within its coplactic graph as
$D$~is'' (because of a similar characterisation for Robinson's bijection,
cf.\ \ref{L-R intro, \Sec4.2}). Alternatively, our bijection~$\pi_1$ could be
taken as basis for the definition of coplactic operations on domino tableaux:
it follows from \ref{L-R intro, theorem~3.3.1} that if $\pi_1(D)=(U,V)$, and a
coplactic operation~$c$ is applicable to~$U$ (and~$V$), then $(c(U),c(V))$ is
again a self-switching tableau pair, whence one could define
$c(D)=\pi_1\inv(c(U),c(V))$. Then the mentioned isomorphism of coplactic
graphs is obvious, but it is not easy to understand directly the effect of a
coplactic operation on a domino tableau. It is the essence of theorem~\pithm\
to establish a link between these two possible definitions of coplactic
operations. As a preliminary step we shall in this section recast the
definition of~\ref{splitting} into a more manageable form that avoids the
choice of a particular reading order.

\subsection Augmented domino tableaux.

As a basis of our definitions we shall use the notion of
$\nu/\kappa$-dominance of domino tableaux, which refines that of Yamanouchi
domino tableaux, and is analogous to the same notion for ordinary tableaux
defined in \ref{L-R intro, 1.4.1}. Definitions using ``companion tableaux'' or
a reading of the entries would, while possible, be less natural in the case of
domino tableaux; nevertheless a simple definition is possible, using
augmentation of the domino tableaux with ordinates added to the entries, as in
the construction of theorem~\yamyamthm.

\proclaim Definition. \augmentdef
An augmentation of a semistandard domino tableau~$D\in\SSDT(\\/\mu)$ for a
skew shape~$\nu/\kappa$ is an assignment of ordinates to all $d\in\Dom(D)$,
such that the conditions below are satisfied; if such an augmentation exists
(which is then unique) $D$ will be called $\nu/\kappa$-dominant. The
combination of the entry~$i$ and ordinate~$j$ of a domino~$d$ will be written
as $i_j$, and $\pos(d)$ will then be denoted by~$\pi(i_j)$.
\defitem Each $i_j$ occurs at most for one domino, and it does so if and only
         if $(i,j)\in Y(\nu/\kappa)$;
\defitem $\pi(i_{j+1})<\pi(i_j)$ whenever
         $\{(i,j),(i,j+1)\}\subset Y(\nu/\kappa)$;
\defitem $\pi((i+1)_j)<\pi(i_j)$ whenever
         $\{(i,j),(i+1,j)\}\subset Y(\nu/\kappa)$.

Conditions \defitemnr1~and~\defitemnr2 imply that among the dominoes with a
fixed entry~$i$, the ordinates increase from right to left by units steps,
starting at~$\kappa_i$; the test for $\nu/\kappa$-dominance then amounts to
verification of condition~\defitemnr3. Like for ordinary tableaux, $D$ can
only be $\nu/\kappa$-dominant if $\wt D=\nu-\kappa$, and we shall say that $D$
is simply $\kappa$-dominant if it is $\nu/\kappa$-dominant for
$\nu=\kappa+\wt{D}$; note that it is still required that $\nu$ be a partition.
We shall show that $D$ is a Yamanouchi domino tableau if and only if it is
$0$-dominant.

The definition above appears to be symmetric with respect to entries~$i$ and
ordinates~$j$, but is not really so due to the requirement that $D$ be a
semistandard domino tableau: ordering dominoes by increasing entries and then
(for equal entries) by decreasing ordinates defines a standard domino tableau.
We shall see that it is possible as well to order dominoes by ordinate first;
these should decrease, and among dominoes with equal ordinates, the entries
should increase. In fact we shall prove the following.

\proclaim Proposition. \specialiseprop
Let $D\in\SSDT(\\/\mu)$ be augmented for~$\nu/\kappa$, and let~`$\prec$' be a
total ordering of the set $\setof i_j: (i,j)\in Y(\nu/\kappa)\endset$ with the
property that $i_j\preceq k_l$ whenever $i\leq k$ and~$j\geq l$. Then there
exists a standard domino tableau~$D_\prec$, called the specialisation of~$D$
for~`$\prec$', with $\Dom(D_\prec)=\Dom(D)$, and in which these dominoes are
added in increasing order for~`$\prec$' of their entry-ordinate combinations
in~$D$.

This means essentially that a domino containing~$i_j$ can be added to the
tableau as soon as the dominoes containing $(i-1)_j$~and~$i_{j+1}$ (if they
exist) are present. The stronger formulation of the proposition is necessary
however: requiring only $(i-1)_j\prec i_j$ and $i_{j+1}\prec i_j$ may fail to
imply $i_j\prec k_l$ for $i<k$ and~$j>l$ by transitivity, as intermediate
entry-ordinate combinations may be absent (by contrast,
definition~\augmentdef\ does imply $\pi(k_l)\leq\pi(i_j)$ whenever $k\geq
i$~and~$l\geq j$, by convexity of skew diagrams). Note that the
orderings~`$\prec$' allowed by the proposition are just the valid reading
orders of the skew diagram~$Y(\nu/\kappa)$. To prove the proposition, we use
the following general lemma (a standard fact of some theory, no doubt),
which allows us to focus on changing the order between just a single pair of
values.

\proclaim Lemma. \bubblesortlemma
Let $X$ be a finite set, $S$ the set (of cardinal $(\Card X)!$) of all total
orderings of~$X$, and $\Gamma$ the graph on~$S$ in which there is an edge
between two orderings if and only if they differ only in the comparison of a
single pair of elements of~$X$ (this is a Cayley graph of $\Sym_{\Card X}$).
Then for any partial ordering~`$<$' on~$X$, the subgraph of\/~$\Gamma$ on the
orderings that extend~`$<$' is connected.

\proof
It is readily checked that any minimal length path in~$\Gamma$ between points
$s,t$ in the subgraph remains within the subgraph. Some may find the following
form of this argument clearer: arranging~$X$ in a linear array according
to~$s$, and then bubble sorting it using~$t$ as sorting criterion, any
comparable pair of elements of~$X$ for~`$<$' remains in the correct relative
order throughout the sorting process.
\QED

\proofof{proposition \specialiseprop}
It will suffice to prove that if the proposition holds for some~`$\prec$', and
$\{i_j,k_l\}$ with $i<k$ and $j<l$ are adjacent for~`$\prec$', then the
proposition holds also for the ordering~`$\prec'$' in which these two
combinations are interchanged. To this end we consider the standard domino
subtableau of~$D_\prec$ consisting of the corresponding two dominoes. Using
edge sequences one easily classifies the possible standard domino tableaux
with the same shape as a given tableau with two dominoes~$d,d'$, in terms of
$|\pos(d)-\pos(d')|$: when this value is~$2$, only one such tableau exists,
when it is~$1$, there are two tableaux, but with different sets of dominoes
(in this case the skew diagram is a $2\times2$ block), and in the remaining
case of a value exceeding~$2$ (since $\pos(d)=\pos(d')$ is impossible), there
are two tableaux, with the same sets of dominoes, which are added in either
order. In view of definition~\augmentdef, we therefore only have to exclude
the possibility that $\pi((i+1)_{j+1})=\pi(i_j)-2$ for some~$i,j$. Assuming
that, one necessarily has $\pi(i_{j+1})=\pi((i+1)_j)=\pi(i_j)-1$; calling this
value~$p$, the fact that the domino containing~$i_{j+1}$ is followed in the
standardisation of~$D$, after some intermediate dominoes, by a domino
containing~$(i+1)_j$ that has the same position~$p$, necessitates the presence
among those intermediate dominoes of a domino with position~$p+2$, and one
with position~$p-2$. Those intermediate dominoes are the ones with entry~$i$
and ordinates~$\leq j$, followed by those with entry~$i+1$ and ordinates~$>j$
(each ordered by decreasing ordinates); therefore the domino with
position~$p+2$ can only be the one containing $i_{j-1}$, while the one with
position~$p-2$ must be the one containing $(i+1)_{j+2}$. But then one has
$\pi((i+1)_j)=\pi(i_{j-1})-2$ and $\pi((i+1)_{j+2})=\pi(i_{j+1})-2$;
taking~$j$ to be either minimal or maximal, one arrives at a contradiction.
\QED

As an instance of this proposition, each augmented domino tableau has a
specialisation for the Kanji order of entry-ordinate pairs ($i_j\prec k_l$
whenever $j>l$, or $j=l$~and~$i<k$). Consequently, the union of a pair of
dominoes $a,b$ respectively containing $i_j$~and~$(i+1)_j$ is a skew diagram;
in particular, $a$ does not extend to the left beyond~$b$ (as in
$\smallsquares\smallsquares{\rtab2 \ribbon0,0:a;0 \ribbon2,1:b;1 }\,$), nor
does~$b$ extend to the right beyond~$a$ (as in
$\smallsquares\smallsquares{\rtab2 \ribbon1,0:a;1 \ribbon2,0:b;0 }\,$).

\proclaim Proposition. \Yamdomprop
A semistandard domino tableau is Yamanouchi if and only if it is $0$-dominant.

\proof
Recall that a domino tableau is Yamanouchi when its column reading is a
Yamanouchi word, and that the column reading orders dominoes by their leftmost
column, and in case of equality by their bottommost row. Yamanouchi words can
be characterised as follows: if one attaches to each letter an ordinate
counting the occurrences of the same letter to its right in the word, then for
any letter~$(i+1)_j$, the letter~$i_j$ occurs to its right. It is easily
verified that if $D$ is a Yamanouchi domino tableau, then attaching these
ordinates of letters to the corresponding entries of~$D$, the conditions of
definition~\augmentdef\ are satisfied: one obtains an augmentation of~$D$ for
$\nu/0$, where~$\nu=\wt D$. Conversely, let $D$ be a semistandard domino
tableau augmented for~$\nu/0$; if~$(i+1)_j$ occurs, then so will~$i_j$, and by
the remarks above, the domino containing~$(i+1)_j$ will precede the one
containing~$i_j$ in the column reading.
\QED

\subsection Coplactic operations on domino tableaux.
\labelsec\coplacdomsec

We now proceed to define coplactic operations on domino tableaux, using
augmentations. A similar description of coplactic operations for ordinary
tableaux can be found in \ref{L-R intro, \Sec3.4}; it essentially consists of
performing a jeu de taquin slide on the set of entry-ordinate pairs, viewed as
coordinates of squares. In the case of domino tableaux we proceed similarly,
but a few additional cases have to be treated. We shall first describe this
procedure, then show that this gives a proper definition of coplactic
operations, and finally show that the operations defined in~\ref{splitting}
produce the same effect as ours.

Like in the case of ordinary tableaux, the coplactic operations
$e_i$~and~$f_i$ only involve the part of the tableau with entries
$i$~and~$i+1$; for the purpose of defining these operations, we may therefore
ignore any other dominoes. So let $D$ be a semistandard domino tableau, all of
whose entries are in~$\{i,i+1\}$, and for the moment consider only the
operation~$e_i$. There exists some~$\kappa$ for which $D$ is $\kappa$-dominant
(it suffices that $\kappa_i-\kappa_{i+1}$ be sufficiently large), and we
equip~$D$ with the augmentation for~$\nu/\kappa$, where~$\nu=\kappa+\wt D$. If
$\kappa_i=\kappa_{i+1}$, then $e_i$ cannot be applied ($D$ is then
$0$-dominant, hence Yamanouchi); otherwise we shall construct a
$\kappa'$-dominant tableau~$E$ with corresponding augmentation, where
$\kappa'$ is obtained from~$\kappa$ by decreasing the part~$\kappa_i$ by~$1$.
It may be that $D$ is already is $\kappa'$-dominant, in which case we shall
get $E=D$ (only ordinates will have changed); if so, the construction is
restarted with the new augmentation found (in other words, $\kappa$ is replaced
by~$\kappa'$). The more interesting case is when $D$ is not $\kappa'$-dominant
(so $\kappa_i-\kappa_{i+1}$ is minimal for $\kappa$-dominance of~$D$); then
obviously $E\neq D$, and we define~$e_i(D)=E$.

We construct~$E$ (with its augmentation) by transforming the augmentation
of~$D$ in a sequence of steps, which as indicated above resemble those of jeu
de taquin, performed on the set of entry-ordinate combinations. At each stage
there is one such combination~$i_j$ with $(i,j)\in Y(\nu/\kappa')$ that does
not occur. Then we consider the combinations $i_{j+1}$ and~$(i+1)_j$; if
neither occurs (so $(i,j)$ is a corner of~$Y(\nu)$), then the process terminates.
Otherwise, the process may proceed by changing one of these combinations
into~$i_j$ (written $i_j\gets i_{j+1}$ or $i_j\gets(i+1)_j$): this simple case
applies when either just one of these two combinations occurs, or when the
union of the dominoes containing these combinations forms a skew diagram; the
change then affects the domino~$d$ of these two for which $\pos(d)$ is larger.

In the remaining case that dominoes exist containing $i_{j+1}$~and~$(i+1)_j$,
but their union is not a skew diagram, one of the following two
transformations is applied, which rearrange several dominoes at once. To allow
a compact display, we subtract~$i$ from all entries and~$j$ from all
ordinates, so that for instance $i_{j+1}$~and~$(i+1)_j$ are drawn respectively
as $0_1$~and~$1_0$. The first possible transformation is
\bigdisplay
  \vertr\: \qquad
  \verconf=
  {\rtab2 \ribbon0,0:0_1;0 \ribbon2,0:1_1;1 \ribbon2,1:1_0;1 }
  \quad\longrightarrow\quad
  {\rtab2 \ribbon1,0:0_1;1 \ribbon1,1:0_0;1 \ribbon2,0:1_0;0 }
  \ .
  \label(\vertrdisp)
$$
The second possible transformation may involve an arbitrarily large number of
dominoes; the largest part of~$D$ that matches the left hand
pattern~$\horconf$ is transformed into the corresponding right hand pattern:
\bigdisplay
  \hortr\: \qquad
  \horconf=
  {\rtab2 \ribbon0,0:0_k;0 \ribbon1,0:1_k;0
         \ribbon0,2:\clap{$\ldots$};0 \ribbon1,2:\clap{$\ldots$};0
	 \ribbon0,4:0_2;0 \ribbon1,4:1_2;0
         \ribbon0,6:0_1;0 \ribbon1,6:1_1;0 \ribbon1,8:1_0;1 }
  \quad\longrightarrow\quad
  {\rtab2 \ribbon1,0:0_k;1
         \ribbon0,1:\hss0_{k-1}\hss;0 \ribbon1,1:\hss1_{k-1}\hss;0
         \ribbon0,3:\clap{$\ldots$};0 \ribbon1,3:\clap{$\ldots$};0
	 \ribbon0,5:0_1;0 \ribbon1,5:1_1;0 
         \ribbon0,7:0_0;0 \ribbon1,7:1_0;0 }\ .
  \label(\hortrdisp)
$$
We define the \newterm{active region} of the application of~$e_i$ in the
simple case as the domino whose entry changes, and otherwise as the pattern
$\verconf$ or~$\horconf$ being transformed.

As an example, consider the following domino tableau~$D$, augmented for
$\\/\mu=(5,5,4)/(2,1,0)$.
$$ D=
{\rtab2
  \ribbon4,1:0_4;1 \ribbon1,5:0_3;0 \ribbon1,7:0_2;1
  \ribbon5,0:1_4;0 \ribbon3,2:1_3;0 \ribbon2,4:1_2;0 \ribbon2,6:1_1;0
  \ribbon4,2:2_3;0 \ribbon4,4:2_2;1 \ribbon4,5:2_1;1 \ribbon3,6:2_0;0
}
$$
Application of~$e_0$ is achieved by the sequence of replacements
$0_1\gets0_2\gets0_3\gets1_3\gets1_4$ (where we ignore the dominoes with
entry~$2$); these only involve the simple case of the definition, since each
of the pairs of dominoes containing $\{0_2,1_1\}$, $\{0_3,1_2\}$, and
$\{0_4,1_3\}$ occupies a skew diagram. Application of~$e_1$ to $D$ however
involves the replacement $1_0\gets1_1$, followed by a transformation~$\vertr$
changing the configuration
$${\rtab2 \ribbon0,0:1_2;0 \ribbon2,0:2_2;1 \ribbon2,1:2_1;1 }
\ttext{into}
 {\rtab2 \ribbon1,0:1_2;1 \ribbon1,1:1_1;1 \ribbon2,0:2_1;0 }\,,
$$
and finally the replacement $2_2\gets2_3$; one therefore finds
\bigdisplay
e_0(D)=
{\rtab2
  \ribbon4,1:0;1 \ribbon1,5:0;0 \ribbon1,7:0;1
  \ribbon5,0:1;0 \ribbon3,2:0;0 \ribbon2,4:1;0 \ribbon2,6:1;0
  \ribbon4,2:2;0 \ribbon4,4:2;1 \ribbon4,5:2;1 \ribbon3,6:2;0
}
\ttext{and}
e_1(D)=
{\rtab2
  \ribbon4,1:0;1 \ribbon1,5:0;0 \ribbon1,7:0;1
  \ribbon5,0:1;0 \ribbon3,2:1;0 \ribbon3,4:1;1 \ribbon3,5:1;1
  \ribbon2,6:1;0
  \ribbon4,2:2;0 \ribbon4,4:2;0 \ribbon3,6:2;0
}\ .
$$
Note that retaining the ordinates found after computing~$e_0(D)$, one would
not get an augmentation of that tableau, which is a consequence of ignoring
the entries~$2$ in the construction, in order to isolate an application
of~$e_0$: while it would in a way be more natural to end the process with
$1_3\gets2_3$ instead of~$1_3\gets1_4$, that would compute $e_1(e_0(D))$ in
one go rather than~$e_0(D)$ (we shall say a bit more about this below). As we
are in fact computing just~$e_i(D)$, there is no real need to continue the
process after the ``hole'' has propagated into row~$i+1$ of~$Y(\nu/\kappa')$,
either after a replacement $i_j\gets (i+1)_j$ or a transformation
$\vertr$~or~$\hortr$; this helps however to relate to the inverse
operation~$f_i$. That operation is defined by reversing the above process in
the obvious ways: entry-ordinate pairs are modified as in an outward jeu de
taquin slide, with in addition the inverse transformations of
$\vertr$~and~$\hortr$. Applying $f_i$ to $e_i(D)$ constitutes a step-by-step
reversal of the application of~$e_i$ to~$D$.

Observe in the computation of~$e_1(D)$ that the augmented tableau contains
$\def\x#1{\hss\scriptstyle\,#1\hss}\smallsquares{\rtab2 \ribbon0,0:\x{1_3};0
\ribbon1,0:\x{2_3};0 \ribbon1,2:\x{2_2};1 }\,$, matching $\horconf$ for
$i_j=1_2$ and $k=1$; nevertheless $\hortr$ does not apply, since the
``hole'' never reaches~$1_2$ (there is no replacement $1_1\gets1_2$ because
the dominoes containing $\{1_2,2_1\}$ in~$D$ do not form a skew diagram).

We shall now proceed to show that the definition above is a proper one, i.e.,
that it handles all cases that can arise, and that the result is always a
semistandard domino tableau. Assume first that the simple case of the
definition applies at all steps; then the result can be seen to be
semistandard by comparison of dominoes, as follows. To facilitate the
formulation, let us define relations `$\lh$'~and~`$\lv$' between dominoes,
where $a\lh b$ and $a\lv b$ both imply that $\pos(a)<\pos(b)$ and that
$a\union b$ is a skew diagram, in addition to which $a\lh b$ means that there
is a standard domino tableau that adds the domino~$a$ and then~$b$, and $a\lv
b$ that there is a standard domino tableau that adds $b$ and then~$a$. Then
either $a\lv b\lh c$ or $a\lh b\lv c$ imply $\pos(c)-\pos(a)\geq 3$, and hence
both $a\lh c$ and $a\lv c$. We may assume there is a replacement of the
form~$i_j\gets(i+1)_j$ (otherwise there is nothing to prove); denote the
domino affected by that replacement by~$d$, and by~$d_{k,l}$ any other domino
containing~$k_l$ \emph{at~that~moment}. Now from the original
contents of the dominoes, and the tests leading to the
replacement~$i_j\gets(i+1)_j$, one deduces $d_{i+1,j+1}\lv d_{i,j+1}\lh d\lh
d_{i+1,j-1}\lv d_{i,j-1}$ (with the proviso that some of the~$d_{k,l}$ may not
exist, but the dominoes actually present form a contiguous subsequence); the
required semistandardness follows.

In the case that at some point the dominoes containing $i_{j+1}$~and~$(i+1)_j$
do not unite to a skew diagram, we argue using the same notation, again
referring to the contents of dominoes at this point. By
proposition~\specialiseprop, there is a standard domino tableau that adds the
dominoes of~$D$ in the following order of their original contents:
$i_{j+1},(i+1)_{j+1},i_j,(i+1)_j$ (where $(i+1)_{j+1}$~and/or~$i_j$ may be
absent). The domino originally containing~$i_j$, if present, is now
$d_{i,j-1}$, and since $d_{i+1,j}\lh d_{i+1,j-1}\lv d_{i,j-1}$, it can be
interchanged in the ordering with~$d_{i+1,j}$; this leaves $d_{i+1,j+1}$ as
the only domino possibly intervening between $d_{i,j+1}$~and~$d_{i+1,j}$, and
since by assumption $d_{i+1,j}$ cannot be added directly after~$d_{i,j+1}$,
the domino $d_{i+1,j+1}$ must be present. It then follows that one either has
a configuration~$\verconf$, or $\horconf$ for~$k=1$; after extending the
latter if possible to a maximal match of~$\horconf$, one finds the
transformation $\vertr$~or~$\hortr$ that is to be applied. The fact that the
result of this transformation is semistandard is verified easily: for~$\vertr$
one uses that presence of~$(i+1)_{j+2}$ implies that of~$i_{j+2}$;
for~$\hortr$, one uses that the presence of~$(i+1)_{j+k+1}$ implies that
of~$i_{j+k+1}$, plus the maximality of the match of~$\horconf$, and
$d_{i+1,j}\lh d_{i,j-1}$ as found above.

For tableaux with entries $i$~and~$i+1$ only, the final entry-ordinate pairs
define a proper augmentation, so $e_i(D)$ is $\kappa'$-dominant, as claimed.
To obtain a similar result in the presence of arbitrary entries, the
construction simulating jeu de taquin can be extended across more rows
of~$Y(\nu/\kappa')$, as was hinted at above. This requires extending the
configurations in~(\vertrdisp) downwards, in a similar way as the
configurations in~(\hortrdisp) are extended leftwards, again choosing a
maximal match for the transformation to be applied. The effect of this
extended construction on the underlying domino tableau (forgetting the
ordinates) will be application of a sequence $e_i,e_{i+1},\ldots$ of
operations with increasing index. If $D$ is $\kappa$-dominant and $\kappa'$ is
obtained from~$\kappa$ by decreasing the part~$\kappa_i$ by~$1$, then the
length of this sequence is the minimal one for which the final result is
$\kappa'$-dominant; this is the same as what happens in the case of ordinary
tableaux \ref{L-R intro, theorem~3.4.1}. Since this construction will not be
used in the sequel, we leave the precise formulations and verifications to the
interested reader; they are quite similar to the ones above.

We now verify that the operations of~\ref{splitting} coincide with ours; we
may assume all entries to be in~$\{i,i+1\}$. First observe that there can be
no other dominoes sharing a column with the active region; therefore it splits
the column reading word into a prefix, an active part, and a suffix. Then the
prefix is a Yamanouchi word, and the suffix is an anti-Yamanouchi (i.e.,
reverse anti-dominant) word; to see this one reasons as in the proof of
proposition~\Yamdomprop, using the ordinates at the time of the
replacement~$i_j\gets(i+1)_j$, or of the transition $\vertr$~or~$\hortr$. For
the simple case this already shows that application of~$e_i$ to~$D$ induces
application of~$e_i$ on the reverse column reading word; for the other cases
it reduces the verification to that for $\verconf$~and~$\horconf$. The
definitions of~\ref{splitting} produce
\bigdisplay
  {\rtab2 \ribbon0,0:0;0 \ribbon2,0:1;1 \ribbon2,1:1;1 }
  \quad\longrightarrow\quad
  {\rtab2 \ribbon0,0:0;0 \ribbon2,0:1;1 \ribbon2,1:0;1 }
  \quad\buildrel\textstyle R_1\over\longrightarrow\quad
  {\rtab2 \ribbon1,0:0;1 \ribbon1,1:0;1 \ribbon2,0:1;0 }
$$
in the first case and
\bigdisplay
 {\rtab2 \ribbon0,0:0;0 \ribbon1,0:1;0
         \ribbon0,2:\clap{$\ldots$};0 \ribbon1,2:\clap{$\ldots$};0
	 \ribbon0,4:0;0 \ribbon1,4:1;0 \ribbon1,6:1;1 }
  \ \longrightarrow\ 
 {\rtab2 \ribbon0,0:0;0 \ribbon1,0:1;0
         \ribbon0,2:\clap{$\ldots$};0 \ribbon1,2:\clap{$\ldots$};0
	 \ribbon0,4:0;0 \ribbon1,4:1;0 \ribbon1,6:0;1 }
  \ \buildrel\textstyle R_2\over\longrightarrow\ 
 {\rtab2 \ribbon0,0:0;0 \ribbon1,0:1;0
         \ribbon0,2:\clap{$\ldots$};0 \ribbon1,2:\clap{$\ldots$};0
	 \ribbon1,4:0;1 \ribbon0,5:0;0 \ribbon1,5:1;0 }
  \ \buildrel\textstyle R_2\over\to\cdots\buildrel\textstyle R_2\over\to\ 
 {\rtab2 \ribbon1,0:0;1
         \ribbon0,1:0;0 \ribbon1,1:1;0
         \ribbon0,3:\clap{$\ldots$};0 \ribbon1,3:\clap{$\ldots$};0
	 \ribbon0,5:0;0 \ribbon1,5:1;0 }
$$
in the second case; in either case the overall effect matches that of
$\vertr$~respectively~$\hortr$.

 \par

\newsection Self-switching tableau pairs and coplactic operations.

In this section we shall prove that under the bijection~$\pi_1$ of
proposition~\selfswitchprop, application of a coplactic operation to a
semistandard domino tableau~$D$ corresponds to application of the same
operation to each component of the associated pair~$(U,V)=\pi_1(D)$. This will
immediately prove theorem~\pithm: the second component in the bijection of
\ref{splitting, Theorem~7.3} is constructed from the sequence of raising
operations that transform~$D$ into a Yamanouchi domino tableau, and since
$U$~and~$V$ have coplactic graphs isomorphic to that of~$D$, the result is the
same as the second component of Robinson's bijection applied to $U$~or~$V$,
which is the Young tableau associated to either of them by jeu de taquin
(cf.~\ref{L-R intro}). This commutation of~$\pi_1$ with coplactic operations
is easy to verify in the basic configurations affected by these operations, in
the various situations that can arise. Nevertheless, the proof in the general
case is difficult, because of the way $\pi_1$ is defined in terms of the
standardisation of~$D$, and the fact that in terms of the standardisation, the
effect of coplactic operations is not localised (although it is in terms of
dominoes and entries). The reduction to the case that $D$ contains entries
$i$~and~$i+1$ is fairly easy, but to further reduce to a situation simple
enough to allow a uniform description of the entire computation of~$\pi_1(D)$
presents difficulty, as relatively distant dominoes intervene in the
standardisation. We shall resolve this problem by showing that other
specialisations of an augmentation of~$D$ can also be used to
compute~$\pi_1(D)$, in particular the one for the Kanji order; the
augmentation will be the one used in defining the coplactic operation.

\subsection Self-switching tableau pairs and augmentation.

In order to show that alternative orders of dominoes can be used to
compute~$\pi_1(D)$, let us consider conditions under which the interchange of
two successive entries in a standard domino tableau leads under application
of~$\pi_0$ to the interchange of the corresponding entries in each component
of the pair of self-switching standard tableaux. Interchanging the order of
two successive entries in a standard domino or ordinary tableau (when this is
possible) will affect exactly one partition in the associated chain; the case
we are interested in is when in the tableau switching family
computing~$\pi_0$, the change of the member~$\\^{[k,k]}$ causes a change of
the members $\\^{[k,j]}$~and~$\\^{[i,k]}$ for all~$i,j$, and of no others. Now
if in a tableau switching family $(\\^{[i,j]})_{i\in I,j\in J}$ with
$I=\{k-1,k,k+1\}$, none of the skew diagrams $Y(\\^{[k+1,j]}/\\^{[k-1,j]})$ is
a domino, then another tableau switching family can be obtained by changing
each partition~$\\^{[k,j]}$ into the unique other partition between
$\\^{[k+1,j]}$~and~$\\^{[k+1,j]}$. This amounts to observing that if during
some sequence of slides applied to a skew tableau of size~$2$, the two squares
never become adjacent, then the same sequence of slides would arise if the
entries of the two squares were interchanged (since apparently the entries are
never being compared); a formal proof can also easily be given by analysing
all possible cases for~$J=\{l,l+1\}$. In order that our tableau switching
family for~$\pi_0$ satisfy the given description, it is necessary that the
subfamily with $i\in\{k-1,k,k+1\}$ be of this type (and its ``transpose''
subfamily with $j\in\{k-1,k,k+1\}$ will then be similar), and that it is still
of this type when its members $\\^{[k-1,k]}$,~$\\^{[k,k]}$, and~$\\^{[k+1,k]}$
are changed (this addition means the dominoes $Y(\\^{[k,k]}/\\^{[k-1,k-1]})$
and $Y(\\^{[k+1,k+1]}/\\^{[k,k]})$ are non-adjacent, so that their entries can
indeed be interchanged). Recall that for ordinary tableaux one has the notions
of $\nu/\kappa$-dominance and augmentations, with the analogue of
proposition~\specialiseprop; calling orderings~`$\prec$' described in that
proposition ``valid'', we have the following.

\proclaim Proposition. \specswitchprop
Let $D$ be a semistandard domino tableau, and~$(U,V)=\pi_1(D)$. Then for any
skew shape~$\nu/\kappa$, each of\/ $U$~and~$V$ is $\nu/\kappa$-dominant if and
only if~$D$ is; if so, for any valid ordering~`$\prec$' of~$Y(\nu/\kappa)$,
the specialisations $D_\prec$,~$U_\prec$, and~$V_\prec$ of their augmentations
for~$\nu/\kappa$ are related by $\pi_0(D_\prec)=(U_\prec,V_\prec)$.

We shall derive the second part first. To be precise, we assume that
$U$~and~$V$ are $\nu/\kappa$-dominant, and attach ordinates to the dominoes
of~$D$ in such a way that conditions \defitemnr1~and~\defitemnr2 of
definition~\augmentdef\ are satisfied; we then show that there is a standard
domino tableau~$D_\prec$ in which the dominoes of~$\Dom(D)$ are added in
increasing order for~`$\prec$', and that $\pi_0(D_\prec)=(U_\prec,V_\prec)$.
This is obvious when `$\prec$' is the Semitic order (the specialisations are
then standardisations), so by lemma~\bubblesortlemma, it suffices to prove
that the stated property for a valid ordering~`$\prec$' implies the same for
another such ordering~`$\prec'$' that differs from it by the exchange of two
consecutive entry-ordinate pairs for~`$\prec$', which are necessarily of the
form $\{i_j,k_l\}$ with $i<k$~and~$j<l$. Since jeu de taquin preserves
$\nu/\kappa$-dominance and is compatible with taking specialisations (cf.\
\ref{pictures, theorem~5.1.1}), the squares containing $i_j$~and~$k_l$ are not
adjacent in any tableau obtained by jeu de taquin from~$U_\prec$ (one has
$\pi(i_j)<\pi(i_l)<\pi(k_l)$, where $\pi(x_y)=\pos(s)$ for the square~$s$
containing~$x_y$). It follows that the tableau switching family for the
computation of~$\pi_0(D_\prec)$ is of the type described above with respect to
the interchange of $i_j$~and~$k_l$, whence
$\pi_0(D_{\prec'})=(U_{\prec'},V_{\prec'})$ as required.

It remains to prove the first part of the proposition. To show that in the
argument above condition~\defitemnr3 of definition~\augmentdef\ is satisfied
for the ordinates added to~$D$, we apply what we have just proved for the case
that `$\prec$' is the Kanji order. Then any pairs $i_j$~and~$(i+1)_j$ are
consecutive, and by the argument used to prove proposition~\selfswitchprop,
the positions of the dominoes containing these pairs in~$D$ are ordered in the
same sense as positions of the corresponding squares of~$U$, which gives the
required condition. For proving conversely that $U$~and~$V$ are
$\nu/\kappa$-dominant whenever~$D$ is, we cannot yet apply the relation
$\pi_0(D_\prec)=(U_\prec,V_\prec)$, so a direct proof is needed, using
standardisations to compute $(U,V)=\pi_1(D)$. To show that the conditions
required for the candidate augmentations of $U$~and~$V$ for $\nu/\kappa$ hold,
we may separately verify them for each restriction to a subtableau with
entries in~$\{i,i+1\}$ only. Since such a restriction is obtained by jeu de
taquin from a component of the tableau switching pair computed from the
similar restriction of~$D$, we may assume that $D$ itself has entries
in~$\{i,i+1\}$ only. Moreover, we may prove $\nu/\kappa$-dominance for any
intermediate tableau in the tableau switching family considered (establishing
$(U,V)=X(U,V)$), in place of $U$~and~$V$; we choose the one at the transition
between the slides into squares with entry~$i$ and with entry~$i+1$. For
any $i_j$~and~$(i+1)_j$, it follows from definition~\augmentdef\ and
proposition~\specialiseprop\ that the outward square of the domino of~$D$
containing~$i_j$ has a position that is strictly larger than that of the
inward square of the domino containing~$(i+1)_j$; from these, the squares of
the indicated intermediate tableau containing the same combinations are
obtained by a sequence of outward respectively inward jeu de taquin slides,
whence their positions are ordered in the same sense.
\QED

If we accept the validity of theorem~\pithm, which will be established in the
next subsection, the proposition allows us to generalise~(\phieq), or more
precisely the equivalent identity
$\<\phi^2(s_\chi),s_\nu>=\eps_2(\chi)\Card{\Yam_2(\chi,\nu)}$:

\proclaim Corollary. \chipsicorr
For all skew shapes~$\chi,\psi$ one has
$$\<\phi^2(s_\chi),s_\psi>=\eps_2(\chi)D(\chi,\psi), \label(\phichieq)
$$
where $D(\chi,\psi)$ is the number of $\psi$-dominant tableaux
in~$\SSDT(\chi)$.

\proof
If the bijection of \ref{splitting, Theorem~7.3} maps $D\in\SSDT(\chi)$ to
$(Y,P)$ with $Y\in\Yam_2(\chi,\\)$ and $P\in\SST(\\)$, then it follows from
proposition~\specswitchprop, theorem~\pithm, and the fact that jeu de taquin
preserves $\psi$-dominance, that $P$ is $\psi$-dominant if and only if $D$ is.
Then $D(\chi,\psi)$ is equal to the sum, taken over all~$\\\in\Part$, of
$\Card{\Yam_2(\chi,\\)}$ times the number of $\psi$-dominant tableaux
in~$\SST(\\)$. The latter number equals $\<s_\\,s_\psi>$, so that we get $
\eps_2(\chi)D(\chi,\psi)
=\sum_{\\\in\Part}\<\phi^2(s_\chi),s_\\>\<s_\\,s_\psi>
=\<\phi^2(s_\chi),s_\psi> $.
\QED

\subsection Commutation of $\pi_1$ with coplactic operations.

Let us state explicitly what we mentioned earlier, that coplactic operations
can be applied to self-switching tableau pairs; this follows by applying
\ref{L-R intro, theorem~3.3.1} twice to  the identity $X(U,V)=(U,V)$.

\proclaim Proposition.
If $c$ is a coplactic operation and $(U,V)$ a self-switching tableau pair,
then $c$ can be applied to~$U$ if and only if it can be applied to~$V$; if so,
$(c(U),c(V))$ is a self-switching tableau pair.
\QED

This subsection is devoted to proving the following theorem; as we have
observed above, it justifies the definition of coplactic operations for
domino tableaux, and it implies theorem~\pithm.

\proclaim Theorem. \complacselfswthm
If $c$ is a coplactic operation, $D$ a semistandard domino tableau, and
$(U,V)=\pi(D)$, then $c$ can be applied to $U$~and~$V$ if and only if it can
be applied to~$D$, and then $(c(U),c(V))=\pi_1(c(D))$.

The proof of this theorem will consist of a reduction in two steps to a
situation sufficiently simple to be verified by an explicit computation. Both
reduction steps are similar in nature, although the second one is much more
subtle; they consist of dividing~$D$ into three regions, of which only the
middle one is affected by the coplactic operation, and showing that $U$
and~$V$ have a similar division. In order to allow the desired conclusion to
be drawn we need a small lemma, whose formulation is more awkward than that
its claim is actually difficult. Recall that standard (domino) tableaux $S,T$,
which are just chains of partitions, can be joined together into~$S|T$ when
their shape match, i.e., the last partition of~$S$ coincides with the first
one of~$T$. We shall use a similar concatenation of semistandard tableaux, but
do not insist that the result correspond to another semistandard tableau
(although it will in the final applications of the construction); it is
defined just as a formal symbol that need only determine a well defined
standard tableau. Then we extend~$\pi_1$ to an operation~$\pi_2$ defined for
such concatenations, using the fact that tableau switching essentially uses
only the structure of standard tableaux. Formally we define
$\pi_2(D_0|\cdots|D_k)=(U_{0,0}|\cdots|U_{0,k}\,,\,V_{0,k}|\cdots|V_{k,k})$,
where $U_{i,j}$~and~$V_{i,j}$ are defined inductively for $0\leq i\leq j\leq
k$ by $(U_{i,i},V_{i,i})=\pi_1(D_i)$, followed by
$(U_{i,j},V_{i,j})=X(V_{i,j-1},U_{i+1,j})$ for~$i<j$. Then we have

\proclaim Lemma. \pitwolemma
If $\pi_2(D_0|D_1|D_2)=(U_0|U_1|U_2\,,\,V_0|V_1|V_2)$, and if $D_1$ satisfies
$\pi_1(c(D_1))=(c(U),c(V))$ where $(U,V)=\pi_1(D_1)$, for a coplactic
operation~$c$, then
$\pi_2(D_0|c(D_1)|D_2)=(U_0|c(U_1)|U_2\,,\,V_0|c(V_1)|V_2)$.

\proof
This follows immediately from the definition, using \ref{L-R intro,
theorem~3.3.1} where appropriate.
\QED

This lemma directly allows the first reduction, which is to the case that $D$
has entries in~$\{i,i+1\}$ only, when $c=e_i$ or~$c=f_i$. The subtableaux
$D_0$,~$D_1$, and~$D_2$ of~$D$ are those whose dominoes have entries~$<i$,
in~$\{i,i+1\}$, and~$>i+1$ respectively; then $(U_0|U_1|U_2\,,\,V_0|V_1|V_2)$
corresponds to $(U,V)=\pi_1(D)$, and clearly $c(U)$ and~$c(V)$ respectively
correspond to $U_0|c(U_1)|U_2$ and~$V_0|c(V_1)|V_2$. Here a concatenation of
semistandard tableaux is said to correspond to another semistandard tableau if
their outer shapes are the same, and the set of dominoes with their entries of
the latter is the (disjoint) union of those of the components of the
concatenation. In this case the relation is particularly simple, since the
concatenation of the standardisations form the standardisation of the combined
tableau, but that will not be the case for the second reduction.

Assume now that $D$ has entries in~$\{i,i+1\}$ only, and that $c=e_i$ (the
case~$c=f_i$ will follow by reversal). The computation of~$e_i(D)$ is achieved
by modifying an augmentation of~$D$ into an augmentation of~$e_i(D)$. We shall
specify below values $k$~and~$l$, which are such that the entry-ordinate
combinations $i_k$ and~$(i+1)_l$ lie on the path of this ``slide'' (i.e., each
is the missing combination at some point during the construction); assume for
the moment that they are given. Let $D_0$,~$D_1$, and~$D_2$ be the
semistandard subtableaux of~$D$ of those dominoes which, at the point in the
construction where $i_k$ is absent, have ordinates~$>l$,
in~$\{l,l-1,\ldots,k\}$, and~$<k$ respectively; similarly let $E_0$,~$E_1$,
and~$E_2$ be the subtableaux of~$e_i(D)$ of those dominoes whose ordinates
satisfy the same conditions respectively, but at the point where $(i+1)_l$ is
absent. Since any changes other than just changing ordinates take place
between these two points in the construction, we have $E_0=D_0$,
$E_1=e_i(D_1)$, and $E_2=D_2$, these are indeed subtableaux by
proposition~\specialiseprop, and $D_0|D_1|D_2$ and $E_0|E_1|E_2$ respectively
correspond to~$D$ and to~$e_i(D)$.

We can now apply lemma~\pitwolemma, but in order for the result to be of any
use, we must ensure that concatenations like $U_0|U_1|U_2$ correspond to the
proper semistandard tableau (in the mentioned case to~$U$ in
$(U,V)=\pi_1(D)$); this is where the specific choices of $k$~and~$l$ become
essential (for badly chosen values, $U_0|U_1|U_2$ might not even correspond to
any semistandard tableau at all). The standard domino tableau~$S$ associated
to~$D_0|D_1|D_2$ is far from being the standardisation of~$D$; however,
proposition~\specswitchprop\ allows any specialisation of the augmentation
of~$D$ to be used to compute (the corresponding specialisations of)
$U$~and~$V$. Now $C$ is almost such a specialisation, but not quite; the only
wrinkle is that, since one must of course use the original augmentation of~$D$
in which the ordinates of all entries~$i$ of~$D_2$ are one more than at the
point where $i_k$ is absent, the domino~$d_2$ of~$D_2$ initially
containing~$i_k$ (if present) is added in~$C$ \emph{after} the domino~$d_1$
of~$D_1$ containing~$(i+1)_k$, whereas $i_k\prec(i+1)_k$ in any valid
order~`$\prec$'. We may assume that $i_k$ is present (otherwise we are already
done), which implies the occurrence of~$(i+1)_{k-1}$ (since $i_k$ changes
into~$i_{k-1}$); if the tableau switching family computing $(S,T)=\pi_0(C)$
can be shown to satisfy the conditions mentioned before
proposition~\specswitchprop, with respect to interchange of $d_1$~and~$d_2$,
then $\pi_2(D_0|D_1|D_2)$ will correspond to~$\pi_1(D)$. Similar
considerations will be needed for $E_0|E_1|E_2$ and~$e_i(D)$; here the
dominoes containing $i_l$~and~$(i+1)_l$ in the augmentation of~$e_i(D)$ must
be exchanged.

We start with considerations independent of the choice of~$k$. Since the
domino~$d_2$ can be added after~$d_1$, the two are not adjacent. We shall now
argue that in the tableau switching family computing $\pi_0(C)$, the
conditions required to the right of the diagonal member~$\\^{[m,m]}$ affected
by interchange of $d_1$~and~$d_2$ are always met; these conditions are that
$Y(\\^{[m+1,n]}/\\^{[m-1,n]})$ is a not domino for any~$n\geq m$. If any of
them were, it would be a horizontal domino, since $\pos(d_1)<\pos(d_2)$. We
claim that this cannot happen, because the square~$s$ of the final difference
$\\^{[m+1,n]}/\\^{[m,n]}$ (i.e., for maximal~$n$) lies in a row strictly above
that of $\\^{[m,m]}/\\^{[m-1,m]}$ (the outer square of~$d_1$); since the
squares of the supposed horizontal domino are obtained from these by inward
respectively outward slides, this gives a contradiction. In fact, we show that
$s$ lies strictly above the outer square~$t$ of the domino~$d_3$ of~$D$
containing~$(i+1)_{k-1}$, which lies weakly above the outer square of~$d_1$.
We augment~$D_2$ as at its definition, so that $d_2$ contains~$i_{k-1}$ and
$d_3$ contains~$(i+1)_{k-1}$; by proposition~\specswitchprop, the second
component of~$\pi_1(D_2)$ has a similar augmentation. Its square
containing~$i_{k-1}$, which is~$s$, then lies strictly above the one
containing~$(i+1)_{k-1}$ which, being obtained by an outward slide from~$t$,
already at the bottom of its column in~$D$, lies in the same row as~$t$.

It remains to show that the conditions required to the left of~$\\^{[m,m]}$
are also met. For this we need an additional assumption, namely that the
square of $\\^{[m+1,m]}/\\^{[m,m]}$ (the inner square of~$d_2$) lies strictly
above any square of~$D_1$; it then lies in particular above any square
obtained by inward slides from that of $\\^{[m,m]}/\\^{[m-1,m]}$, which
suffices as above. To verify the assumption, we choose~$k$ to be the largest
ordinate for which the inner square of the domino containing~$i_k$ in the
initial augmentation of~$D$ lies in a row strictly above the active region for
the application of~$e_i$, or if no such~$k$ exists, the one for which the
combination $i_k$ is initially absent. The considerations for~$e_i(D)$ are
symmetric, and lead to choosing for~$l$ the smallest ordinate for which the
outer square of the domino containing~$(i+1)_l$ in the final augmentation
of~$e_i(D)$ lies in a row strictly below the active region, or if no such~$l$
exists, the one for which the combination $(i+1)_l$ is finally absent. With
these choices, we have achieved a reduction of the verification of
theorem~\complacselfswthm\ for~$D$ to that for~$D_1$.

\proofof{theorem~\complacselfswthm}
By the above reductions we may assume that $D$ has entries in~$\{i,i+1\}$
only, and that all dominoes lie in the rows containing the active region; we
may moreover assume that if transforming the augmentation of $D$ into that
of~$e_i(D)$ starts with $i_k$ absent and ends with $(i+1)_l$ absent, then
$k\leq j\leq l$ for all ordinates~$j$. If the active region for the
application of~$e_i$ to~$D$ consists either of a single domino or of a
configuration~$\horconf$, its top right corner contains~$i+1$ in~$D$, and its
bottom left corner contains~$i$ in~$e_i(D)$, so that there can be no
entries~$i$ to the right nor entries~$i+1$ to the left; it then follows that
$D$ is reduced to the active region. In these cases the theorem is verified by
a simple computation, displayed below. Like before we distinguish the
components of the tableau switching pairs by italic and bold entries, and
subtract~$i$ from all entries, so that they become $0$~or~$1$; we also display
ordinates, running from~$0$ up to~$r$, with $r'=r-1$. Horizontal arrows denote
applications of~$e_i$ (to each component in case of tableau switching pairs),
and vertical arrows applications of~$\pi_1$:
\bigdisplay
\def\v{\noalign{\medskip}\downarrow& &\downarrow\cr\noalign{\medskip}}
\displaylines{
 \matrix
 { {\rtab2 \ribbon0,0:1_0;0 } &\to & {\rtab2 \ribbon0,0:0_0;0 }\cr
 \v
  \Young({\it1}_0,{\bf1}_0) &\to & \Young(\it0_0,{\bf0}_0)\cr
 }
\hfil;\hfil
 \matrix
 { {\rtab2 \ribbon1,0:1_0;1 } &\to & {\rtab2 \ribbon1,0:0_0;1 }\cr
 \v
  \Young({\it1}_0|{\bf1}_0) &\to & \Young({\it0}_0|{\bf0}_0)\cr
 }
\hfil;\hfil
\matrix
{
  {\rtab2 \ribbon0,0:0_r;0 \ribbon1,0:1_r;0
	 \ribbon0,2:0_{r'}\hss;0 \ribbon1,2:1_{r'}\hss;0
         \ribbon0,4:\clap{$\ldots$};0 \ribbon1,4:\clap{$\ldots$};0
         \ribbon0,6:0_1;0 \ribbon1,6:1_1;0 \ribbon1,8:1_0;1 }
 &\to
 &{\rtab2 \ribbon1,0:0_r;1
         \ribbon0,1:0_{r'}\hss;0 \ribbon1,1:1_{r'}\hss;0
         \ribbon0,3:\clap{$\ldots$};0 \ribbon1,3:\clap{$\ldots$};0
	 \ribbon0,5:0_1;0 \ribbon1,5:1_1;0 
         \ribbon0,7:0_0;0 \ribbon1,7:1_0;0 }
\cr
\v
  {\rtab1
   \ribbon0,0:{\it0}_r; \ribbon0,3:{\it0}_1;
   \ribbon0,4:{\it1}_0; \ribbon0,5:{\bf0}_r; \ribbon0,8:{\bf0}_1; 
   \ribbon1,0:{\it1}_r; \ribbon1,3:{\it1}_1;
   \ribbon1,4:{\bf1}_r; \ribbon1,7:{\bf1}_1; \ribbon1,8:{\bf1}_0;
   \ribbonlength=2
   \ribbon0,1:\hss\cdots\hss;0 \ribbon0,6:\hss\cdots\hss;0
   \ribbon1,1:\hss\cdots\hss;0 \ribbon1,5:\hss\cdots\hss;0
  }
 &\to
&  {\rtab1
   \ribbon0,0:{\it0}_r; \ribbon0,1:{\it0}_{r'}\hss; \ribbon0,4:{\it0}_0;
   \ribbon0,5:{\bf0}_{r'}\hss; \ribbon0,8:{\bf0}_0; 
   \ribbon1,0:{\it1}_{r'}\hss; \ribbon1,3:{\it1}_0;
   \ribbon1,4:{\bf0}_r; \ribbon1,5:{\bf1}_{r'}\hss; \ribbon1,8:{\bf1}_0;
   \ribbonlength=2
   \ribbon0,2:\hss\cdots\hss;0 \ribbon0,6:\hss\cdots\hss;0
   \ribbon1,1:\hss\cdots\hss;0 \ribbon1,6:\hss\cdots\hss;0
  }
\cr}\,.
\cr
}
$$
The case involving a transformation~$\vertr$ is more complicated. In the top
row of the active region, there may be horizontal dominoes with
entry~$i$ to its right, and in the bottom row there may be dominoes with
entry~$i+1$ to its left; these are accompanied by an equal number of entries
$i+1$ respectively~$i$ in the remaining rows, each of which may be divided
into a number of horizontal and vertical dominoes. Hence this case is governed
by four free parameters; in the diagram below these are determined by the
values of the ordinates $0\leq n<p\leq q\leq r$ (as before we subtract a
common offset from all ordinates to make the smallest one equal to~$0$). We
indicate certain sequences of horizontally repeated dominoes or squares with
equal entries by an ellipsis, but to avoid excessive width, we alternatively
indicate such repetition by displaying just one domino or square, with an
asterisk attached to its entry; in that case the ordinate given is the
leftmost (largest) one of the sequence. Primes attached to ordinates signify
subtraction of~$1$.
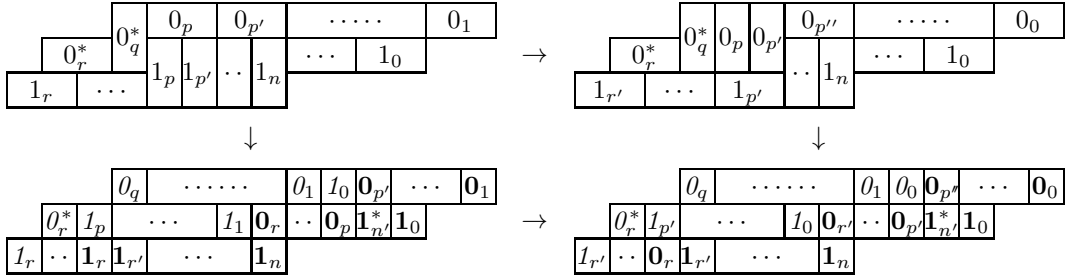
\bigdisplay
\def\v{\noalign{\medskip}\downarrow& &\downarrow\cr\noalign{\medskip}}
 \setbox\strutbox=\hbox
  {\vrule height 1.1\ht\strutbox depth 1.1\dp\strutbox width 0pt }
  \def\[#1]{\hss#1\hss}
 \matrix
 {
 {\rtab2
  \ribbon1,1:0_r^*;0 \ribbon1,3:0_q^*;1
  \ribbon0,4:0_p;0 \ribbon0,6:0_{p'}\hss;0 \ribbon0,12:0_1;0
  \ribbon2,0:1_r;0 \ribbon2,2:\hss\cdots\hss;0 \ribbon2,4:1_p;1
  \ribbon2,5:1_{p'}\hss;1 \ribbon2,6:\cdotp\cdotp;1 \ribbon2,7:1_n;1
  \ribbon1,8:\hss\cdots\hss ;0 \ribbon1,10:1_0;0
  \ribbonlength=4 \ribbon0,8:\hss\cdots\cdotp\cdotp\hss;000
 }
& \to
&
 {\rtab2
  \ribbon1,1:0_r^*;0 \ribbon1,3:0_q^*;1
  \ribbon1,4:0_p;1 \ribbon1,5:0_{p'}\!;1 \ribbon0,6:0_{p''}\hss;0
  \ribbon0,12:0_0;0
  \ribbon2,0:1_{r'}\hss;0 \ribbon2,2:\hss\cdots\hss;0 \ribbon2,4:1_{p'}\hss;0
  \ribbon2,6:\cdotp\cdotp;1 \ribbon2,7:1_n;1
  \ribbon1,8:\hss\cdots\hss ;0 \ribbon1,10:1_0;0
  \ribbonlength=4 \ribbon0,8:\hss\cdots\cdotp\cdotp\hss;000
 }
\cr
\v
 {\rtab1
  \ribbon1,1:{\it0}_r^*; \ribbon0,3:{\it0}_q;
  \ribbon0,8:{\it0}_1; \ribbon0,9:{\it1}_0;
  \ribbon2,0:{\it1}_r; \ribbon2,1:\cdotp\cdotp;
  \ribbon1,2:{\it1}_p; \ribbon1,6:{\it1}_1;
  \ribbon1,7:{\bf0}_r; \ribbon1,8:\cdotp\cdotp; \ribbon1,9:{\bf0}_p;
  \ribbon0,10:{\bf0}_{p'}\!; \ribbon0,13:{\bf0}_1;
  \ribbon2,2:{\bf1}_r; \ribbon2,3:{\bf1}_{r'}\hss; \ribbon2,7:{\bf1}_n;
  \ribbon1,10:{\bf1}_{n'}^*\hss; \ribbon1,11:{\bf1}_0;
  \ribbonlength=2 \ribbon0,11:\cdots;0
  \ribbonlength=3 \ribbon1,3:\cdots;00 \ribbon2,4:\cdots;00
  \ribbonlength=4 \ribbon0,4:\hss\cdots\cdots\hss;000
 }
& \to
&
 {\rtab1
  \ribbon1,1:{\it0}_r^*; \ribbon0,3:{\it0}_q;
  \ribbon0,8:{\it0}_1; \ribbon0,9:{\it0}_0;
  \ribbon2,0:{\it1}_{r'}; \ribbon2,1:\cdotp\cdotp;
  \ribbon1,2:{\it1}_{p'}; \ribbon1,6:{\it1}_0;
  \ribbon1,7:{\bf0}_{r'}\!; \ribbon1,8:\cdotp\cdotp; \ribbon1,9:{\bf0}_{p'}\!;
  \ribbon0,10:{\bf0}_{p^{\prime\!\prime}}\!; \ribbon0,13:{\bf0}_0;
  \ribbon2,2:{\bf0}_r; \ribbon2,3:{\bf1}_{r'}\hss; \ribbon2,7:{\bf1}_n;
  \ribbon1,10:{\bf1}_{n'}^*\hss; \ribbon1,11:{\bf1}_0;
  \ribbonlength=2 \ribbon0,11:\cdots;0
  \ribbonlength=3 \ribbon1,3:\cdots;00 \ribbon2,4:\cdots;00
  \ribbonlength=4 \ribbon0,4:\hss\cdots\cdots\hss;000
 }
\cr
}
$$
The shape of the tableaux is indicated only schematically, and away from the
active region, vertical alignment is mostly accidental; however, in the bottom
left diagram the alignment of ${\it1}_p$~and~${\bf1}_r$ and that of
${\it1}_0$~and~${\bf0}_p$ are always as indicated, as are the corresponding
alignments in the bottom right diagram. It may be verified that the diagram
correctly describes the computation for all possible parameter values.
\QED

 \par

\newsection Moving chains and coplactic operations.

In this section we shall consider the same kind of questions as in the
previous one, but replacing the correspondence~$\pi_1$ by the operation of
moving chains. In will be necessary to limit the type of chains in order to
obtain the desired properties; in particular commutation with coplactic
operations is obtained for the operation of moving open chains.

\subsection Moving chains in augmented domino tableaux.

Let a semistandard domino tableau~$D$ be given, and a chain~$C$
for~$s\in\{s_0,s_1\}$ in~$D$ that is not a forbidden chain, then \ref{edge
sequences, proposition~4.3.1} tells us that it can be moved in~$D$, resulting
in another semistandard domino tableau~$D'$. We suppose that $D$ is
$\nu/\kappa$-dominant, and ask under which circumstances the same is true
for~$D'$; this generalises the question when moving a chain preserves the
Yamanouchi property. It is easy to see that if we augment~$D$
for~$\nu/\kappa$, and move the ordinates of the dominoes in~$C$ along with
that domino, then conditions \defitemnr1~and~\defitemnr2 of
definition~\augmentdef\ are satisfied (in fact the forbidden chains are
precisely those that would violate condition~\defitemnr2); therefore our
question is equivalent to asking whether condition~\defitemnr3 is satisfied.
Fix~$i,j$ for which there are dominoes in~$\Dom(D)$ containing respectively
$i_j$ and $(i+1)_j$, at least one of which lies in~$C$. The positions of their
fixed squares for~$s$ are congruent modulo~$2$, and a violation of
condition~\defitemnr3 can only occur if these positions are equal. If so, we
have a configuration
$$
   {\rtab2 \ribbon0,0:i_j;0 \ribbon1,0:\hss{i^+}_j\hss;0 }\,,
  \label(\blockedconf)
$$
(where $i^+$ denotes~$i+1$), with the fixed squares for~$s$ on its main
diagonal. If the top domino lies in~$C$, then either the bottom domino is its
successor in~$C$, or there is another configuration~(\blockedconf) directly to
the left of the one considered, but with ordinates $j+1$ instead of~$j$;
similarly if the bottom domino lies in~$C$, then either the top domino is its
successor in~$C$, or there is another configuration~(\blockedconf), with
ordinates $j-1$, to the right of this one. It follows that
condition~\defitemnr3 of definition~\augmentdef\ is violated if and only
if~$C$ meets a configuration~(\blockedconf), in which case it is a closed
chain, all of whose dominoes occur in such a configuration. The following
definition covers both such pairs of dominoes and the ones forming a forbidden
chain.

\proclaim Definition. \blockeddef
Let $D$ be an semistandard domino tableau, provided with an augmentation
for~$\nu/\kappa$. A pair of dominoes of~$D$ is called blocked for
$s\in\{s_0,s_1\}$ if the positions of their fixed squares for~$s$ are equal,
and their entry-ordinate combinations correspond to adjacent squares
of~$Y(\nu/\kappa)$. A chain in~$D$ for~$s$ is called blocked if any of its
dominoes occurs in a blocked pair, and unblocked otherwise.

\proclaim Proposition. \augmoveprop
Let $D$ be a semistandard domino tableau, provided with an augmentation
for~$\nu/\kappa$, and $C$ a chain in~$D$ for $s\in\{s_0,s_1\}$. A
$\nu/\kappa$-dominant tableau can be obtained from~$D$ by moving~$C$ if and
only if $C$ is unblocked, which is the case as soon as any domino of~$C$ is not
part of a blocked pair.
\QED

The notions of blocked pairs and blocked chains depend on the augmentation
of~$D$. However, open chains are never blocked, so in particular moving open
chains in preserves the Yamanouchi property; this proves our second main
theorem,~\blowcorethm. Also, whether or not a chain in an augmented tableau is
blocked does not change when other (unblocked) chains are moved, which proves
lemma~\closchindlem.

\proclaim Proposition.
If $D$ is semistandard domino tableau, provided with an augmentation, then any
unblocked chain~$C$ in~$D$ for~$s$ is also a chain in any specialisation of
the augmentation of~$D$.

\proof
Let~$D'$ be the result of moving~$C$ in~$D$, which by
proposition~\augmoveprop\ can be augmented for the same skew shape as~$D$ was;
for every occurring $i_j$, the dominoes of $D$~and~$D'$ containing that
combination have their fixed square for~$s$ in common. Then for any valid
order~`$\prec$' the fact that $D_\prec$~and~$D'_\prec$ are standard domino
tableaux (proposition~\specialiseprop) implies that $C$ is a chain
of~$D_\prec$. We may also formulate this as follows. To determine the
successor and predecessor relations of the chain~$C$, the entries in the
standardisation of~$D$ of certain pairs of dominoes in~$\Dom(D)$ had to be
compared, which is equivalent to comparing their entry-ordinate combinations
using the Semitic order. The corresponding two dominoes in~$\Dom(D')$ are
always adjacent, and since $D'_\prec$ is a well defined standard domino
tableau, this implies that the comparison of their entry-ordinate combinations
using~`$\prec$' gives the same result as using the Semitic order.
\QED

The fact that, by definition~\blockeddef, the entry-ordinate
combinations in a blocked pair of a domino tableau augmented for~$\nu/\kappa$
define a domino inside~$Y(\nu/\kappa)$, is the key to the construction of
theorem~\yamyamthm.

\proofof{theorem~\yamyamthm}
Recall that, for each individual domino~$d$ of an augmented Yamanouchi domino
tableau in~$\Yam_2(\\,\mu)$, we constructed a pair of dominoes forming a
$2\times2$ block inside the diagram~$Y(\dublpart\mu)$, with top-left corner
at~$(2i,2j)$ if $d$ contains~$i_j$. These dominoes, considered in isolation,
can be made into a blocked pair for~$s_1$, by giving them entry-ordinate
combinations corresponding to the coordinates of the squares of~$d$ (using the
inward square of~$d$ for the domino with the larger position). Note that these
entries agree with the ones we assigned in the construction, and that the
assignment of ordinates was also already shown in our example in \Sec2. It is
easy to see that this gives is a bijective correspondence between on one hand
tilings~$L$ of $Y(\\)$ with dominoes, each provided with an entry and
ordinate, satisfying only condition~\defitemnr1 of definition~\augmentdef\ for
$\nu/\kappa=\mu/0$ (without requirement of belonging to a semistandard domino
tableau), and on the other hand tilings~$M$ of~$Y(\dublpart\mu)$ with blocked
pairs of dominoes, with the entries and ordinates again satisfying only
condition~\defitemnr1 of definition~\augmentdef, but for $\nu/\kappa=\\/0$.
What remains to be shown is that $L$ defines an augmentation of a semistandard
domino tableau if and only if $M$ does so.

First note the semistandardness of an assignment of entries to a tiling of a
diagram with dominoes can be tested by comparing adjacent dominoes: when two
dominoes touch along a vertical edge, their entries should increase weakly
from left to right, and when they touch along a horizontal edge, their entries
should increase strictly from top to bottom. This test can be seen to succeed
for $L$ or~$M$ if the other satisfies definition~\augmentdef. For~$M$, we need
only consider edges between different $2\times2$ blocks. For horizontal edges,
the condition follows because both squares of the domino containing~$(i+1)_j$
in~$L$ lie strictly below those of the domino containing~$i_j$. For vertical
edges one similarly uses weak increase of row numbers from~$i_j$ to~$i_{j+1}$,
which fails only in one case when $i_j$~and~$i_{j+1}$ themselves fill a
$2\times2$ block in~$L$, but that case is readily checked to cause no problem
either. Semistandardness of~$L$ is verified similarly, but here row numbers
are divided by~$2$ (and truncated to an integer), which causes strict increase
down columns to be the subtler of the two conditions: if square~$(r,c)$
belongs to a different domino of~$L$ than~$(r+1,c)$, then the domino of~$M$
containing~$r_c$ cannot lie entirely in an even row~$2i$, which forces its row
number to remain strictly smaller than that of the domino
containing~$(r+1)_c$, even after division by~$2$ and truncation. The
verification of conditions \defitemnr2~and~\defitemnr3 of
definition~\augmentdef\ proceed almost by reversal of these arguments, but in
addition to monotonicity of entries one uses strict decrease of ordinates
along rows and weak decrease down columns (which is easily checked), so that
the difference $j-i$ for combinations~$i_j$ decreases strictly both along rows
and down columns (not counting the places where one stays within a domino, of
course). For~$M$ this suffices to establish the mentioned conditions, since
the necessary comparisons of the positions of dominoes can be done based only
on the $2\times2$ blocks containing them (which may be assumed distinct);
for~$L$ one observes that for these comparisons one may compare the contents
of the inward dominoes of the corresponding $2\times2$ blocks in~$M$, or
equivalently the contents associated to the top-left squares of these blocks.
\QED

The fact that $L$~and~$M$ have partition shapes plays no r\^ole at all in the
proof above; these conditions are simply remnants of the fact that we had not
yet introduced the notion of $\nu/\kappa$-dominance when stating
theorem~\yamyamthm, but only the notion of Yamanouchi domino tableau.
We can therefore generalise as follows:

\proclaim Corollary. \domdomcorr
For any pair of skew shapes $\chi,\psi$ with $|\chi|=2|\psi|$, there is a
bijection between the set of $\psi$-dominant domino tableaux~$L$ of shape
$\chi$, and the set of $\chi$-dominant domino tableaux~$M$ of shape
$\dublpart\psi$ for which there is no chain for~$s_1$ that can be moved
without destroying its $\chi$-dominance. The bijection is given by the
construction of theorem~\yamyamthm, applied to the appropriate augmentations;
in particular $2\Spin(L)=|\psi|-\Spin(M)$ (i.e., $M$ has twice as many
horizontal dominoes as $L$ has vertical dominoes).
\QED

In the particular case that $\chi=\\/\mu$ is a horizontal strip (its columns
have length at most~$1$), every domino tableau~$M$ of weight~$\\-\mu$ is
$\chi$-dominant. Then any $M$ in the corollary consists only of forbidden
chains for~$s_1$, and the set of all such~$M$ is in bijection, not just with
the indicated set of domino tableaux~$L$, but more directly (by shrinking each
forbidden chain to a single square) with the set of $M'\in\SST(\psi)$ with $\wt
M'={1\over2}(\\-\mu)$; it is the latter correspondence that was used in the
derivation of equation~(\spingfcanceleq).

\subsection Coplactic operations and moving open chains.

We have now proved all the facts announced in \Sec2, but one more question
deserves to be considered, namely whether coplactic operations commute with
moving open chains. We already know that moving open chains does not affect
the applicability of coplactic operations (since it preserves
$\nu/\kappa$-dominance), and the close analogy to \ref{L-R intro,
theorem~3.3.1} is a strong indication for expecting commutation. That one does
in fact have commutation, as we shall now demonstrate, is therefore not
surprising, but a proof still requires a considerable amount of work, and the
verification of quite a few particular cases.

\proclaim Theorem. \complacmocthm
Coplactic operations on domino tableaux commute with moving of open chains, in
the following sense. If $D\in\SSDT(\chi)$, $C$ is an open chain in~$D$
for~$s\in\{s_0,s_1\}$, and~$D'$ is obtained by moving~$C$ in~$D$, then $c(D)$
is defined if and only if $c(D')$~is, for any coplactic operation~$c$; if so,
the result of moving the open chain~$C'$ of~$c(D)$ starting in the same square
as~$C$ is equal to~$c(D')$.

\proof
We may assume that $c=e_i$, and (since the statement about applicability of
coplactic operations is already proved) that $e_i(D)$ is defined. We set
$E=e_i(D)$, and denote the result of moving~$C'$ in~$E$ by~$E'$; we then have
to show that $E'=e_i(D')$. We provide $D$ with the augmentation that is
transformed into one of~$E$ by the definition of coplactic operations, and
provide $D'$, $E'$, and~$e_i(D')$ with corresponding augmentations.

Assume first that $C=C'$ (i.e., these chains consist of the same sequence of
dominoes; their entries and ordinates in $D$ respectively~$E$ may differ); we
must show that in this case the steps performed in the application of~$e_i$
to~$D$ (replacements of entry-ordinate pairs and possibly a transformation
$\vertr$~or~$\hortr$) are matched exactly in the application of~$e_i$ to~$D'$.
It is clear that if a transformation $\vertr$~or~$\hortr$ occurs, the affected
configuration $\verconf$~or~$\horconf$ is disjoint from~$C$, since all
dominoes of~$C$ must remain in~$E$ in order to have $C=C'$. It therefore
suffices to show that every replacement $i_j\gets i_{j+1}$ or
$i_j\gets(i+1)_j$ during the computation of~$e_i(D)$ has its counterpart in
the computation of~$e_i(D')$. We may assume that both candidate combinations
$i_{j+1},(i+1)_j$ for this replacement are present in the augmentation of~$D$
(and of~$D'$), for otherwise the statement is obvious for lack of choice; let
$\{x,y\}\subset\Dom(D)$ be the dominoes containing these combinations, with
the one of~$x$ actually being replaced (so $\pos(x)>\pos(y)$). Denoting by
$x',y'\in\Dom(D')$ the corresponding dominoes after moving~$C$, it follows
from~$C'=C$ that $x',y'\in\Dom(E')$, and that~$x'$ contains~$i_j$ in the
augmentation of~$E'$, while~$y'$ contains the same combination
($(i+1)_j$~or~$i_{j+1}$) as it does in the augmentation of~$D'$. This implies
that the union of $x'$~and~$y'$ is a skew diagram, and that
$\pos(x')>\pos(y')$, which proves the required statement.

The remainder of the proof, which shows that in all cases with $C\neq C'$ one
has commutation as well, will appeal mainly to the kind of people who like to
study the solutions of chess problems (but those who prefer to find the
solution themselves may stop reading and do just that). We shall start by
listing a number of variants in which the theorem can be seen to hold, and
then show that with $C\neq C'$, there is no way to avoid running into one of
them. From any variant~$\phi$, three others can be obtained by the following
symmetries related to those of the ``rectangle'' of tableaux
$\smash{\bigl[{D\atop D'}~{E\atop E'}\bigr]}$; these will not be listed
separately. A variant denoted by~$\phi^\updownarrow$ is obtained by
interchanging $D\leftrightarrow D'$ and $E\leftrightarrow E'$ and reversing
the movement of the chains $C$~and~$C'$. A variant denoted
by~$\phi^\leftrightarrow$ is obtained by interchanging $D\leftrightarrow E$
and $D'\leftrightarrow E'$ while rotating all tableaux a half-turn and
renumbering entries and ordinates in the opposite order, so as to interchange
the r\^oles of $i$~and~$i+1$. The third variant is $\phi^\times
=\phi^{\updownarrow\leftrightarrow}=\phi^{\leftrightarrow\updownarrow}$.

We do not make the assumption that $D$ has entries in~$\{i,i+1\}$ only, since
suppression of the remaining entries might change the structure of the
chain~$C$. However, we shall focus on a specific part of~$D$ only, in which
all entries are in~$\{i,i+1\}$, and which in all cases has the property that
the dominoes not lying in the explicitly indicated part of~$C$ are member of a
blocked pair for~$s$, whence they cannot be part of~$C$, regardless of the
suppressed context (the situation for $E$, $D'$ and~$E'$ is similar). Some
cases involve transformations $\hortr$, and correspondingly have one or more
parameters describing how often certain elements of the configurations are
repeated horizontally. We include the case of zero repetitions, whenever this
makes sense; we then also consider a replacement $i_j\gets(i+1)_j$ in a
vertical domino as an instance of the transformation $\hortr$ with~$k=0$. We
display entries with their ordinates (as before subtracting offsets to make
the smallest ones~$0$), and indicate the relevant part of the open chain in
each tableau by arrows. The fact that the chains and the effect of~$e_i$ are
as displayed, implies that any repeatable element is included in the display
as often as the context allows. We start giving variants that we shall call
$\alpha(k,l)$, for parameters $k\geq l\geq0$ (variant $\alpha(0,0)$ has
$C=C'$, and may be excluded).
\bigdisplay
\def\tweek#1{\mskip-0.5\thinmuskip#1\mskip-0.5\thinmuskip}
 \matrix
{ D\:\
 {\rtab2
  \ribbon1,0:0_k;0 \ribbon2,0:1_k;0
  \ribbon1,2:\ldots;0 \ribbon2,2:\ldots ;0
  \ribbon1,4:\hss0_{l+1}\hss;0 \ribbon2,4:\hss1_{l+1}\hss;0
  \ribbon2,6:1_l;1
  \ribbon0,7:0_l;0 \ribbon1,7:\hss1_{l\tweek-1}\!\hss;0
  \ribbon0,9:\ldots;0 \ribbon1,9:\ldots ;0
  \ribbon0,11:0_1;0 \ribbon1,11:1_0;0
  \arrow1,6:r \arrow1,8:r \arrow1,10:r \arrow1,12:u
  \arrow0,11:l \arrow0,9:l \arrow0,7:l
 }
& \buildrel \textstyle e_i \over \longrightarrow
& E\:\
 {\rtab2
  \ribbon2,0:\raise2pt\hbox{$0_k$};1
  \ribbon1,1:\hss0_{k\!-\!1\!}\!\hss;0 \ribbon2,1:\hss1_{k\!-\!1}\hss;0
  \ribbon1,3:\ldots;0 \ribbon2,3:\ldots ;0
  \ribbon1,5:0_l;0 \ribbon2,5:1_l;0
  \ribbon0,7:\hss0_{l-1}\hss;0 \ribbon1,7:\hss1_{l-1}\hss;0
  \ribbon0,9:\ldots;0 \ribbon1,9:\ldots ;0
  \ribbon0,11:0_1;0 \ribbon1,11:1_0;0
  \arrow2,5:l \arrow2,3:l \arrow2,1:l
  \arrow1,0:r \arrow1,2:r \arrow1,4:r \arrow1,6:u
 }
\cr
\noalign{\smallskip}
\Big\downarrow\rlap{$C$} & \frame{\frame{$\alpha(k,l)$}} &
\Big\downarrow\rlap{$C'$}\cr
\noalign{\smallskip}
 D'\:\
 {\rtab2
  \ribbon1,0:0_k;0 \ribbon2,0:1_k;0
  \ribbon1,2:\ldots;0 \ribbon2,2:\ldots ;0
  \ribbon1,4:\hss0_{l+1}\hss;0 \ribbon2,4:\hss1_{l+1}\hss;0
  \ribbon0,6:0_l;0 \ribbon1,6:1_l;0
  \ribbon0,8:\ldots;0 \ribbon1,8:\ldots ;0
  \ribbon0,10:0_1;0 \ribbon1,10:1_1;0
  \ribbon1,12:\,1_0\!;1
  \arrow1,6:d \arrow1,8:l \arrow1,10:l \arrow1,12:l
  \arrow0,11:r \arrow0,9:r \arrow0,7:r
 }
&\buildrel \textstyle e_i \over \longrightarrow
&E'\:\
 {\rtab2
  \ribbon1,0:0_k;0 \ribbon2,0:\hss1_{k\!-\!1\!}\!\hss;0
  \ribbon1,2:\ldots;0 \ribbon2,2:\ldots ;0
  \ribbon1,4:\hss0_{l\tweek+1}\hss;0 \ribbon2,4:\hss1_l\hss;0
  \ribbon1,6:0_l;1
  \ribbon0,7:\hss0_{l-1}\hss;0 \ribbon1,7:\hss1_{l-1}\hss;0
  \ribbon0,9:\ldots;0 \ribbon1,9:\ldots ;0
  \ribbon0,11:0_1;0 \ribbon1,11:1_0;0
  \arrow2,5:r \arrow2,3:r \arrow2,1:r
  \llap{\arrow1,0:d ~}\arrow1,2:l \arrow1,4:l \arrow1,6:l
 }
\cr
}
$$
Note that variant $\alpha(k,l)^\times$ coincides with $\alpha(k,k-l)$.
Then one has variants $\beta(k)$ for $k>0$, obtained by modifying
$\alpha(k,0)$ slightly: upward arrows at the end of chains $C$~and~$C'$
in~$\alpha(k,0)$ are replaced by rightward arrows in~$\beta(k)$, which does
not affect the essential characteristics of the configuration.
\bigdisplay
 \matrix
{ 
 {\rtab2  \global\jmax=7
  \ribbon0,0:0_k;0 \ribbon1,0:1_k;0
  \ribbon0,2:\ldots;0 \ribbon1,2:\ldots ;0
  \ribbon0,4:\hss0_1;0 \ribbon1,4:1_1;0
  \ribbon1,6:1_0;1
  \arrow0,6:r
 }
& \buildrel \textstyle e_i \over \longrightarrow
& 
 {\rtab2  \global\jmax=7
  \ribbon1,0:\raise2pt\hbox{$0_k$};1
  \ribbon0,1:\hss0_{k\!-\!1\!}\!\hss;0 \ribbon1,1:\hss1_{k\!-\!1}\hss;0
  \ribbon0,3:\ldots;0 \ribbon1,3:\ldots ;0
  \ribbon0,5:0_0;0 \ribbon1,5:1_0;0
  \arrow1,5:l \arrow1,3:l \arrow1,1:l
  \arrow0,0:r \arrow0,2:r \arrow0,4:r \arrow0,6:r
 }
\cr
\noalign{\smallskip}
\Big\downarrow\rlap{$C$} & \frame{\frame{$\beta(k)$}} &
\Big\downarrow\rlap{$C'$} \cr
\noalign{\smallskip}
 {\rtab2
  \ribbon0,0:0_k;0 \ribbon1,0:1_k;0
  \ribbon0,2:\ldots;0 \ribbon1,2:\ldots ;0
  \ribbon0,4:0_1;0 \ribbon1,4:1_1;0
  \ribbon0,6:1_0;0
  \arrow0,6:d
 }
&\buildrel \textstyle e_i \over \longrightarrow
&
 {\rtab2
  \ribbon0,0:0_k;0 \ribbon1,0:\hss1_{k\!-\!1\!}\!\hss;0
  \ribbon0,2:\ldots;0 \ribbon1,2:\ldots ;0
  \ribbon0,4:0_1;0 \ribbon1,4:1_0;0
  \ribbon0,6:0_0;0
  \arrow1,5:r \arrow1,3:r \arrow1,1:r
  \llap{\arrow0,0:d ~}\arrow0,2:l \arrow0,4:l \arrow0,6:l
 }
\cr
}
$$
A third family of variants is $\gamma(k)$ for $k\geq1$. Here both a
transformation~$\hortr$ (for $e_i$ applied to~$D$) and a transformation
$\vertr$ (for $e_i$ applied to~$D'$) occur. Similar variants~$\gamma'(k)$ are
obtained from~$\gamma(k)$ by changing the direction at the beginning of the
chains $C$~and~$C'$.
\bigdisplay
\dimen0=\ht\strutbox \advance\dimen0\dp\strutbox
 \matrix
{ \noalign{\kern-.5\dimen0}
 {\rtab2
  \ribbon2,0:0_k;0 \ribbon3,0:1_k;0
  \ribbon2,2:\ldots;0 \ribbon3,2:\ldots ;0
  \ribbon3,4:1_1;1
  \ribbon1,5:0_1;1 \ribbon3,5:1_0;1
  \arrow1,5:l
 }
& \buildrel \textstyle e_i \over \longrightarrow
& 
 {\rtab2
  \ribbon1,5:0_0;1   \ribbon3,5:1_0;1
  \ribbon3,0:\raise2pt\hbox{$0_k$};1
  \ribbon2,1:\ldots;0 \ribbon3,1:\ldots;0
  \ribbon2,3:0_1;0 \ribbon3,3:1_1;0
  \arrow3,3:l \arrow3,1:l \arrow2,0:r \arrow2,2:r
  \arrow1,5:d \arrow3,5:l \arrow2,4:u
 }
\cr
\noalign{\smallskip}
\Big\downarrow\rlap{$C$} & \clap{\frame{\frame{$\gamma(k)$}}} &
\Big\downarrow\rlap{$C'$} \cr
\noalign{\smallskip\kern-.5\dimen0}
 {\rtab2
  \ribbon1,0:0_k;0 \ribbon2,0:1_k;0
  \ribbon1,2:\ldots;0 \ribbon2,2:\ldots ;0
  \ribbon2,4:1_1;1
  \ribbon0,4:0_1;0   \ribbon2,5:1_0;1
  \arrow0,5:u
 }
&\buildrel \textstyle e_i \over \longrightarrow
& 
 {\rtab2
  \ribbon1,4:0_1;1
  \ribbon1,5:0_0;1   \ribbon2,4:1_0;0
  \ribbon1,0:0_k;0 \ribbon2,0:\hss1_{k\!-\!1\!}\!\hss;0
  \ribbon1,2:\ldots;0 \ribbon2,2:\ldots ;0
  \arrow2,3:r \arrow2,1:r \llap{\arrow1,0:d ~}\arrow1,2:l \arrow1,4:l
  \arrow0,5:u \arrow2,5:u
 }
\cr
}
\kern1cm
 \matrix
{ \noalign{\kern.5\dimen0}
 {\rtab2
  \ribbon1,0:0_k;0 \ribbon2,0:1_k;0
  \ribbon1,2:\ldots;0 \ribbon2,2:\ldots ;0
  \ribbon2,4:1_1;1
  \ribbon0,5:0_1;0   \ribbon2,5:1_0;1
  \arrow0,5:l
 }
& \kern-.5\dimen0\buildrel \textstyle e_i \over \longrightarrow
& 
 {\rtab2
  \ribbon0,5:0_0;0   \ribbon2,5:1_0;1
  \ribbon2,0:\raise2pt\hbox{$0_k$};1
  \ribbon1,1:\ldots;0 \ribbon2,1:\ldots;0
  \ribbon1,3:0_1;0 \ribbon2,3:1_1;0
  \arrow2,3:l \arrow2,1:l \arrow1,0:r \arrow1,2:r
  \arrow0,5:d \arrow2,5:l \arrow1,4:u
 }
\cr
\noalign{\smallskip}
\Big\downarrow\rlap{$C$} & \clap{\frame{\frame{$\gamma'(k)$}}} &
\Big\downarrow\rlap{$C'$} \cr
\noalign{\smallskip\kern-.5\dimen0}
 {\rtab2 \global\jmax=6
  \ribbon1,0:0_k;0 \ribbon2,0:1_k;0
  \ribbon1,2:\ldots;0 \ribbon2,2:\ldots ;0
  \ribbon2,4:1_1;1
  \ribbon0,4:0_1;0   \ribbon2,5:1_0;1
  \arrow0,5:r
 }
& \kern-.5\dimen0\buildrel \textstyle e_i \over \longrightarrow
& 
 {\rtab2 \global\jmax=6
  \ribbon1,4:0_1;1
  \ribbon1,5:0_0;1   \ribbon2,4:1_0;0
  \ribbon1,0:0_k;0 \ribbon2,0:\hss1_{k\!-\!1\!}\!\hss;0
  \ribbon1,2:\ldots;0 \ribbon2,2:\ldots ;0
  \arrow2,3:r \arrow2,1:r \llap{\arrow1,0:d ~}\arrow1,2:l \arrow1,4:l
  \arrow0,5:r \arrow2,5:u
 }
\cr
}
$$
The final variants involve transformations~$\vertr$ in applying~$e_i$ both
to~$D$ and to~$D'$, but at different (overlapping) places; they have no
parameters. The direction of $C$~and~$C'$ can be varied at both ends, giving
$4$ such variants: $\delta$,~$\delta'$ and~$\delta''$ displayed below,
and~$\delta'^\times$.
\bigdisplay
\dimen0=\ht\strutbox \advance\dimen0\dp\strutbox
\def\squareframe#1{\clap{\frame{\frame{\hbox to\dimen0{\strut\hss$#1$\hss}}}}}
 \matrix
{ 
 {\rtab2
  \ribbon1,2:0_1;1 \ribbon1,0:0_2;0 \ribbon3,0:1_2;1
  \ribbon3,1:1_1;1 \ribbon3,2:1_0;1 
  \arrow1,2:l \arrow1,0:d \arrow3,0:d
 }
& \!\!\buildrel \textstyle e_i \over \longrightarrow \!\!
& 
 {\rtab2
  \ribbon1,2:0_0;1 \ribbon3,2:1_0;1 \ribbon3,0:1_1;0
  \ribbon2,0:0_2;1 \ribbon2,1:0_1;1 
  \arrow1,2:d \arrow3,2:l \arrow3,0:d
 }
\cr
\noalign{\smallskip}
\Big\downarrow\rlap{$C$} & \squareframe\delta &
\Big\downarrow\rlap{$C'$} \cr
\noalign{\smallskip}
 {\rtab2
  \ribbon0,1:0_1;0 \ribbon1,0:0_2;1 \ribbon3,0:1_2;1
  \ribbon2,1:1_1;1 \ribbon2,2:1_0;1 
  \arrow0,2:u \arrow0,0:r \arrow2,0:u
 }
& \!\!\buildrel \textstyle e_i \over \longrightarrow \!\!
& 
 {\rtab2
  \ribbon1,2:0_0;1 \ribbon2,1:1_0;0 \ribbon3,0:1_1;1
  \ribbon1,0:0_2;1 \ribbon1,1:0_1;1 
  \arrow0,2:u \arrow2,2:u \arrow2,0:r
 }
\cr
}
\qquad
 \matrix
{
 {\rtab2
  \ribbon1,2:0_1;0 \ribbon1,0:0_2;0 \ribbon3,0:1_2;1
  \ribbon3,1:1_1;1 \ribbon3,2:1_0;1 
  \arrow1,2:l \arrow1,0:d \arrow3,0:d
 }
& \kern-.5\dimen0 \buildrel \textstyle e_i \over \longrightarrow \!\!
& 
 {\rtab2
  \ribbon1,2:0_0;0 \ribbon3,2:1_0;1 \ribbon3,0:1_1;0
  \ribbon2,0:0_2;1 \ribbon2,1:0_1;1 
  \arrow1,2:d \arrow3,2:l \arrow3,0:d
 }
\cr
\noalign{\smallskip}
\Big\downarrow\rlap{$C$} & \squareframe{\delta'} &
\Big\downarrow\rlap{$C'$} \cr
\noalign{\smallskip}
 {\rtab2 \global\jmax=3 
  \ribbon0,1:0_1;0 \ribbon1,0:0_2;1 \ribbon3,0:1_2;1
  \ribbon2,1:1_1;1 \ribbon2,2:1_0;1
  \arrow0,2:r \arrow0,0:r \arrow2,0:u
 }
& \kern-.5\dimen0 \buildrel \textstyle e_i \over \longrightarrow \!\!
& 
 {\rtab2 \global\jmax=3
  \ribbon1,2:0_0;1 \ribbon2,1:1_0;0 \ribbon3,0:1_1;1
  \ribbon1,0:0_2;1 \ribbon1,1:0_1;1
  \arrow0,2:r \arrow2,2:u \arrow2,0:r
 }
\cr
}
\quad
 \matrix
{ 
 {\rtab2
  \ribbon0,3:0_1;0 \ribbon0,1:0_2;0 \ribbon2,1:1_2;1
  \ribbon2,2:1_1;1 \ribbon2,3:1_0;1
  \arrow0,3:l \arrow0,1:d \arrow2,1:l
 }
&\kern-.5\dimen0 \buildrel \textstyle e_i \over \longrightarrow
 \kern-.5\dimen0
& 
 {\rtab2
  \ribbon0,3:0_0;0 \ribbon2,3:1_0;1 \ribbon2,1:1_1;0
  \ribbon1,1:0_2;1 \ribbon1,2:0_1;1
  \arrow0,3:d \arrow2,3:l \arrow2,1:l
 }
\cr
\noalign{\smallskip}
\Big\downarrow\rlap{$C$} & \squareframe{\delta''} &
\Big\downarrow\rlap{$C'$} \cr
\noalign{\smallskip}
 {\rtab2 \global\jmax=4
  \ribbon0,2:0_1;0 \ribbon1,1:0_2;1 \ribbon2,0:1_2;0
  \ribbon2,2:1_1;1 \ribbon2,3:1_0;1
  \arrow0,3:r \arrow0,1:r \arrow2,1:u
 }
&\kern-.5\dimen0 \buildrel \textstyle e_i \over \longrightarrow
 \kern-.5\dimen0
& 
 {\rtab2 \global\jmax=4
  \ribbon1,3:0_0;1 \ribbon2,2:1_0;0 \ribbon2,0:1_1;0
  \ribbon1,1:0_2;1 \ribbon1,2:0_1;1
  \arrow0,3:r \arrow2,3:u \arrow2,1:r
 }
\cr
}
$$

We now proceed to show that $C'\neq C$ implies that, up to the mentioned
symmetries, one of these variants arises. We label the cases that can arise in
the application of~$e_i$ to~$D$ as follows: case~$a$ means a replacement
$i\gets(i+1)_j$ occurs in a horizontal domino; case~$b$ means that a
transformation~$\vertr$ occurs; case~$c$ means that a transformation~$\hortr$
occurs (which we take to include a replacement $i\gets(i+1)_j$ in a vertical
domino). In all cases the active region for the application of~$e_i$ is
rectangular, and of its top-left and bottom-right corner squares exactly one
is fixed for~$s$; by applying the symmetry $\phi\to\phi^\leftrightarrow$ if
necessary, we assume that it is the top left square.

Assuming $C\neq C'$ as we do, there are two possibilities for the first point
of divergence between $C$~and~$C'$: \ (1)~some domino of~$C$ does not occur
in~$\Dom(E)$ (due to rearrangement of dominoes), and \ (2)~the successor~$y$
in~$C$ of some domino~$x$ does occur in~$\Dom(E)$, but it is not the successor
of~$x$ in~$C'$. Recall from \ref{edge sequences, 4.1} that successor of~$x$ is
determined by comparing its entry to that of the domino occupying its
``discriminant square'', which is one step in the anti-diagonal direction away
from the fixed square (for~$s$) of~$x$ (if the discriminant square lies
outside~$D$, then the test depends only on the shape of~$D$, and we cannot
have possibility~(2)). Therefore, possibility~(2) requires that the entry
either of~$x$ itself or of the domino containing the discriminant square is
changed by~$e_i$, and in such a way that the result of the comparison changes.
Now in each of the cases there is a unique fixed square~$t_0$ for~$s$ for
which application of~$e_i$ changes the entry of its domino: $t_0$~is the left
square of the active domino in case~$a$, in case~$b$ it is the right square of
the middle row of~$\verconf$, and in case~$c$, it is the top-right corner
of~$\horconf$. Let $t_0$ be part of a domino~$x_0\in\Dom(D)$
containing~$(i+1)_j$, and let $t_1$ be the discriminant square of~$x_0$ (it
lies to the bottom-left of~$t_0$ in case~$a$, and to the top-right of~$t_0$
otherwise); if in~$D$ the square $t_1$ is part of a domino containing
either~$(i+1)_{j+1}$ (case~$a$) or~$i_j$ (other cases), then we define~$x_1$
to be that domino. If possibility~(2) applies, the crucial comparison must be
between the entries of $x_0$~and~$x_1$; moreover, in case~$c$ it must be that
$x=x_1$ (since $x_0\notin\Dom(E)$), but in case~$b$, even if~$x_1$ is defined,
its successor lies in~$\verconf$ and therefore not in~$\Dom(E)$, which
excludes possibility~(2) altogether.

Suppose first that we are in case~$a$; then $\Dom(D)=\Dom(E)$, so we have
possibility~(2). Since $x_0$ undergoes a replacement~$i_j\gets(i+1)_j$ and
$x_1$ contains~$(i+1)_{j+1}$ in~$D$, there is a domino~$x_2$
containing~$i_{j+1}$ in~$D$, with $x_1\lv x_2\lh x_0$. Since $x_0$ is
horizontal, the only possible configuration is
$\smallsquares\def\x#1{\hss\raise1pt\hbox{$x_#1$}\hss}\smash{\rtab2
\ribbon1,0:\x1;0 \ribbon0,0:\x2;0 \ribbon0,2:\x0;0 }$, so $x_1$~and~$x_2$
form a blocked pair for~$s$, which forces~$x=x_0$. By including any further
blocked pairs for~$s$ with entries $i$~and~$i+1$ occurring directly to the
left of $x_1$~and~$x_2$, we find an occurrence of~$\beta(k)^\updownarrow$.

Suppose next that we are in case~$b$. Now the bottom left domino of~$\verconf$
(containing~$(i+1)_{j+1}$) is the successor for~$s$ of the top domino
(containing~$i_{j+1}$). It may be that the remaining domino~$x_0$
of~$\verconf$ (containing~$(i+1)_j$) forms a blocked pair for~$s$ with a
domino containing $(i+1)_{j-1}$ in~$D$; if it does, $x_1$ is defined, and one
gets one of the variants $\delta$,~$\delta'$, $\delta'^\times$, or~$\delta''$.
Otherwise $x_1$ may still be defined, but then since its contents is replaced
$i_{j-1}\gets i_j$ when applying~$e_i$, it forms a $2\times2$ block of
horizontal dominoes with a domino containing $(i+1)_{j-1}$ in~$D$. Reasoning
similarly we may find more such $2\times2$ blocks with entries $i$~and~$i+1$
directly to the right, but the chain through~$x_0$ eventually turns back and
passes through the top and bottom left domino of~$\verconf$; if $k$ blocks
were found, we have a variant $\gamma(k+1)^\times$ or~$\gamma'(k+1)^\times$.

Suppose finally that we are in case~$c$. If $x_1$ is not defined we have
possibility~(1), and the domino where $C$~and~$C'$ diverge must be~$x_0$
(being in~$\horconf$ and not part of a blocked pair); this leads to a variant
$\alpha(k,0)$~or~$\beta(k)$. Otherwise, since $e_i$ causes a
replacement~$i_{j-1}\gets i_j$ in~$x_1$, there is a domino~$x_2$ containing
$(i+1)_{j-1}$ in~$D$, with $x_0\lh x_2\lv x_1$, and it is the successor
of~$x_0$ for~$s$ in~$D$. If vertical, $x_2$ forms a blocked pair for~$s$
with~$x_0$, so that $C$ does not meet~$\horconf$, and we must have
possibility~(2); this leads to variant $\gamma(k)$~or~$\gamma'(k)$. If $x_2$
is horizontal, then so is~$x_1$, and the two form a $2\times2$ block. As above
we may need to include more such blocks with entries $i$~and~$i+1$ to the
right, but the chain through~$x_0$ eventually turns back to pass
through~$x_1$, excluding possibility~(2); this leads to a
variant~$\alpha(k,l)$ with~$l>0$.
\QED

For the record we mention that more grandiose versions of the commutation
theorems \complacselfswthm\ and~\complacmocthm\ can be obtained, by replacing
the coplactic operations by the extended construction simulating jeu de taquin
on all entry-ordinate pairs that was indicated in subsection~\coplacdomsec.
The argument that reduces these commutation theorems to the ones given is the
same as used for \ref{L-R intro, corollary~3.4.3}: the number of coplactic
operations that correspond to a single jeu de taquin slide can be read off
from the coplactic graph, and it therefore does not change when $\pi_1$ is
applied or an open chain is moved. The theorems also have a consequence that
does not refer to coplactic operations at all; we do not know a direct proof
of it.

\proclaim Corollary.
If $D'$ is obtained from~$D$ by moving open chains, then $\pi(D)=\pi(D')$.

\proof
By theorems \complacmocthm, \complacselfswthm, and
\ref{L-R intro, theorem~3.3.1}, we can reduce to the case that $D$ is a
Yamanouchi domino tableaux, and $D'$ will be one as well. But then
$\pi(D)=\Can{\wt D}=\Can{\wt D'}=\pi(D')$.
\QED

%
 \par

\newsection Concluding remarks.

The complements we have given to the facts stated in~\ref{splitting} result in
a fairly complete and elegant theory for domino tableaux, which forms a nice
parallel to the theory for ordinary tableaux, as described in~\ref{L-R intro}.
Yet, one may well wonder why all this works out so nicely, and if there is
some deeper algebraic significance behind the combinatorial facts. The most
general and natural algebraic interpretation of (augmented) domino tableaux
that we know of, stated in corollary~\chipsicorr, seems to cry for a
generalisation to higher plethysm operators, involving \r-ribbon tableaux
for~$r>2$. However, none of the constructions of this paper generalise well to
\r-ribbon tableaux, and no reasonable generalisation of the definition of
augmented (or even just Yamanouchi) domino tableaux seems to give the correct
values (e.g., it is hard to decide which tableaux should correspond to
$\<\phi^3(s_{(4,4,4)}),s_{(2,2)}>=1=\<\phi^3(s_{(3,3,3,3)}),s_{(2,2)}>$). The
interested reader may check that every analysis in this paper depends on
properties of domino tableaux that fail for \r-ribbon tableaux with~$r>2$ (or,
as in the case of~$\pi_1$, on constructions that have no analogue at all). One
may surmise that the special place of domino tableaux has some relation with
classical groups of types other than~$A_n$, rather than just with plethysms;
such a connection certainly does exist (see~\ref{Garfinkle}, where even
moving chains play a r\^ole), but at present we cannot indicate any concrete
links.

\Supplement References.

\input \jobname.ref

\iftoc
  \closeout\tocfile \vfil\eject
  \def\tocitem#1=#2\onpage#3.
    {\par\line{\hbox to 1.5pc{\bf #1\hss}\quad#2\dotfill#3}\ignorespaces}
  \input \jobname.toc

\fi

\bye